%% file: bare_jrnl_JETSPE.tex
\documentclass[journal,onecolumn,draftclsnofoot,]{IEEEtran}
\usepackage{color, colortbl}
\usepackage{cite}
\usepackage{mathtools}
\usepackage{amsmath,bm}
\usepackage{amsfonts,bbold}
\usepackage{caption}
\usepackage[T1]{fontenc}
\usepackage[utf8]{inputenc}
\usepackage{subfigure}
\usepackage{tikz}
\usepackage{amsthm}
\usepackage{pgfplots}
\usepackage{wrapfig}
\usepackage[export]{adjustbox}
\usepackage{siunitx}
\usepackage[american,cuteinductors,smartlabels]{circuitikz}
\usetikzlibrary{calc}
\ctikzset{bipoles/thickness=1}
\newtheorem{remark}{Remark}
\definecolor{LightCyan}{rgb}{0.88,1,1}
\newcommand{\tb}{\color{black}}

\newcommand\blue[1]{{\color{black}#1}}

\DeclareMathOperator{\diff}{d}
\DeclareMathOperator{\sat}{sat}
\DeclareMathOperator{\sgn}{sgn}
\DeclareMathOperator{\If}{if}

\providecommand{\norm}[1]{\lVert#1\rVert}

\definecolor{light-gray}{gray}{0.96}
\begin{document}
\title{Interactions of Grid-Forming Power Converters and Synchronous Machines}

\author{Ali~Tayyebi,~
        Dominic~Gro\ss,~
        Adolfo~Anta,~
        Friederich~Kupzog~and Florian~D\"{o}rfler

\thanks{This work was partially funded by the independent research fund of the
    the power system digitalization group at the Electric Energy Systems (EES)
    competence unit of the Austrian Institute for Technology (AIT), ETH Zürich funds, and by the European Unions Horizon 2020 research and innovation programme under grant agreement No. 691800. This article reflects only the authors views and the European Commission is not responsible for any use that may be made of the information it contains. A. Tayyebi is with AIT, 1210 Vienna, Austria, and also with the Automatic Control Laboratory, ETH Z\"{u}rich, Switzerland. A. Anta and F. Kupzog are with AIT. D. Gro\ss, and F. D\"{o}rfler are with the Automatic Control Laboratory, ETH Z\"{u}rich, 8092 Z\"{u}rich, Switzerland; Email: \{ali.tayyebi-khameneh,adolfo.anta,friederich.kupzog\}@ait.ac.at, \{grodo,dorfler\}@ethz.ch.}%

\thanks{Corresponding author's contact information: Ali Tayyebi, PhD Candidate at AIT and ETH Z\"{u}rich, Giefinggasse 2, 1210 Vienna, Austria, Phone: +43 6648251413, Fax: +43 505506390, Email: ali.tayyebi-khameneh@ait.ac.at}}
\maketitle 
\vspace{-1cm}
\renewcommand{\abstractname}{\vspace{-\baselineskip}}
\begin{abstract}
\noindent\textbf{\emph{Abstract}---}An inevitable consequence of the global power system transition towards nearly 100\% renewable-based generation is the loss of conventional bulk generation by synchronous machines, their inertia, and accompanying frequency and voltage control mechanisms. This gradual transformation of the power system to a low-inertia system leads to critical challenges in maintaining system stability. Novel control techniques for converters, so-called grid-forming strategies, are expected to address these challenges and replicate functionalities that so far have been provided by synchronous machines. We present a low-inertia high-fidelity case study that includes synchronous machines and models of grid-forming converters. We study interactions between synchronous machines and converters and analyze the response of various grid-forming control approaches to contingencies, i.e., large changes in load and the loss of a synchronous machine. Our case study highlights the positive impact of the grid-forming converters on frequency stability and analyze the potential limitations of each control technique when interacting with synchronous machines.  Our studies also analyze how and when the interaction between the fast grid-forming converter, the dc source \blue{and ac} current limitations, and the slow synchronous machine dynamics contributes to system instability. \blue{Lastly, we introduce an effective solution to address the instability issues due to the GFCs ac and dc current limitation.}  
\end{abstract}

\section{Introduction}
At the heart of the energy transition is the change in generation technology; from fossil fuel based thermal generation to converter interfaced renewable generation \cite{milano_foundations_2018}. One of the major consequences of this transition towards a nearly  100\% renewable system is the gradual loss of synchronous machines (SMs), their inertia, and control mechanisms. This loss of the rotational inertia changes the nature of the power system to a low-inertia network resulting in critical stability challenges \cite{milano_foundations_2018,SP98,kundur1994power}. On the other hand, low-inertia power systems are characterized by large-scale integration of generation interfaced by power converters, allowing frequency and voltage regulation at much faster time-scales compared to SMs \cite{tayyebi_grid-forming_2018,milano_foundations_2018}.

Indeed, power converters are already starting to provide new ancillary services, modifying their active and reactive power output based on local measurements of frequency and voltage. However, because of the dependency on frequency measurements these \emph{grid-following} control techniques only replicate the instantaneous inertial response of SMs after a contingency with a delay and result in degraded performance on the time scales of interest \cite{noauthor_high_2017}. To resolve this issue, \emph{grid-forming converters} (GFCs) are envisioned to be the cornerstone of future power systems. Based on the properties and functions of SMs, it is expected that grid-forming converters must support load-sharing/drooping, black-start, inertial response, and hierarchical frequency/voltage regulation. While these services might not be necessary in a future converter-based grid, a long transition phase is expected, where SMs and GFCs must be able to interact and ensure system stability.

Several grid-forming control strategies have been proposed in recent years \cite{tayyebi_grid-forming_2018}. \emph{Droop control} mimics the speed droop mechanism present in SMs and is a widely accepted baseline solution \cite{chandorkar_control_1993}. As a natural further step, the emulation of SM dynamics and control led to so-called \emph{virtual synchronous machine} (VSM) strategies \cite{darco_virtual_2013,zhong_synchronverters:_2011,synchronverter2014}. Recently, \emph{matching} control strategies that exploit structural similarities of converters and synchronous machine and \emph{match} their dynamic behavior have been proposed \cite{cvetkovic_modeling_2015,arghir_grid-forming_2017,curi_control_2017,Catalin2019}. In contrast, virtual oscillator control (VOC) uses GFCs to mimic the synchronizing behavior of Li\'enard-type oscillators and can globally synchronize a converter-based power system~\cite{sinha_uncovering_2017}. However, the nominal power injection of VOC cannot be specified. This limitation is overcome by \emph{dispatchable virtual oscillator control} (dVOC)  \cite{colombino_global_2017,gros_effect_2018,SG-SI-BJ-MC-DG-FD:18} that ensures synchronization to a pre-specified operating point that satisfies the ac power flow equations.%

In this article, we explicitly consider the dynamics of the converter dc-link capacitor, the response time of the dc power source, and its current limits. \blue{We review four different grid-forming control strategies and combine them with standard low-level cascaded control design accounting for the ac voltage control and the ac current limitation and control \cite{PPG07}}. We compare various performance aspects of GFC control techniques in an electromagnetic transients (EMT) simulation of the IEEE 9-bus test system, namely: 1) their impact on the frequency performance metrics e.g., nadir and rate of change of frequency (RoCoF) \cite{PM18,PSF2017,entsoefreq_2018,poolla_placement_2018}, 2) their interaction with high-fidelty SM models, 3) their response under large load disturbance magnitudes, \blue{4) their behavior when imposing dc and ac current limitations,} and 5)  their response to the loss of SM. Furthermore, we provide an insightful qualitative analysis of the simulation results. The models used in this work are available online \cite{model}.

This comparative study highlights the positive impact of GFCs on improving standard power system frequency stability metrics. \blue{Moreover, we observe that limiting the GFCs dc or ac current accompanied by the interaction of fast converters and slow synchronous machine dynamics can destabilize VSMs, droop control, and dVOC, while matching control appears to be unaffected.} 
Furthermore, we reveal a potentially destabilizing interaction between the fast synchronization of GFCs and the slow response of SMs (see \cite{gros_effect_2018,Uros2019} for a similar observation on the interaction of GFCs and line dynamics). Lastly, this study shows that an all-GFCs (i.e., no-inertia) system can exhibit more resilience than a mixed SM and GFC (i.e., low-inertia) system with respect to the large load variations.  

The remainder of this article is structured as follows: Section \ref{modeling} reviews the modeling approach. Section \ref{GFCC} presents the system dynamics and adopted grid-forming control techniques. The simulation-based analysis of our comparative study is discussed in Section \ref{simulation}. Section \ref{analysis} presents a qualitative analysis of the observations made in case studies. Our concluding statements and agenda of future work are reported in Section \ref{conclusion}. And the choice of control parameters is described in Appendix \ref{app:tunning}.

\section{Model Description}\label{modeling}
Throughout this study, we use a test system comprised of power converters and synchronous machines. This section describes the models of the individual devices and components \cite{model}.
{\tb
\subsection{Converter Model}\label{sec:ConverterModel}
\begin{figure}[b!]
    \centering    
    \includegraphics[clip, trim=0.25cm 0cm 0cm 0cm,width=0.8\textwidth]{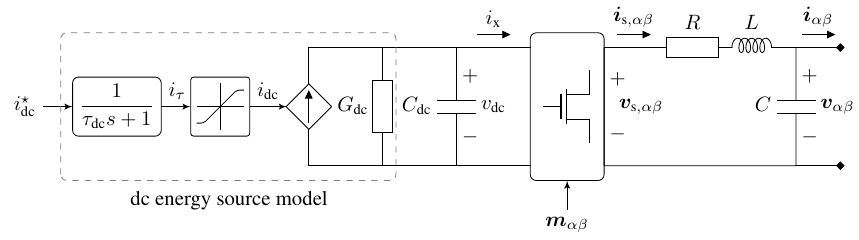}
    \caption{Converter model in $\alpha\beta$-coordinates with detailed dc energy source model based on \eqref{eqs:model}-\eqref{eq:source_sat}.\label{fig:ConverterModel}}
\end{figure}
To begin with, we consider the converter model illustrated in Figure \ref{fig:ConverterModel} in $\alpha\beta$-coordinates \cite{YI10,Catalin2019} 
\begin{subequations}\label{eqs:model}
    \begin{align}
    C_{{\text{\upshape{dc}}}}{\dot{v}}_{{\text{\upshape{dc}}}}&=i_{\text{\upshape{dc}}}{\tb -G_\text{\upshape{dc}}v_\text{\upshape{dc}}}-i_{\text{\upshape{x}}},\label{eqs:model1}\\
    L{\dot{\bm{i}}_{\tb\text{\upshape{s}},\alpha\beta}}&=\bm{v}_{\tb\text{\upshape{s}},\alpha\beta}-R\bm{i}_{\tb\text{\upshape{s}},\alpha\beta}-\bm{v}_{\tb\alpha\beta},\\
    C{\dot{\bm{v}}_{\tb\alpha\beta}}&=\bm{i}_{\tb\text{\upshape{s}},\alpha\beta}-\bm{i}_{\tb\alpha\beta},
    \end{align}
\end{subequations}
where $C_{{\text{\upshape{dc}}}}$ denotes the dc-link capacitance, $G_{{\text{\upshape{dc}}}}$ models dc losses, and, $L$, $C$, and $R$ respectively denote the filter inductance, capacitance, and resistance. Moreover, $v_{\text{\upshape{dc}}}$ represents the dc voltage, $i_{\text{\upshape{dc}}}$ is the current flowing out of the controllable dc current source, $\bm{m}_{\tb\alpha\beta}$ denotes the modulation signal of the full-bridge averaged switching stage model, $i_{\text{\upshape{x}}} \coloneqq (1/2)\bm{m}^\top_{\tb\alpha\beta} \bm{i}_{\tb\text{\upshape{s}},\alpha\beta}$ denotes the net dc current delivered to the switching stage, and $\bm{i}_{\tb\text{\upshape{s}},\alpha\beta}$ and $\bm{v}_{\tb\text{\upshape{s}},\alpha\beta}\coloneqq(1/2) \bm{m}_{\tb\alpha\beta} v_{{\text{\upshape{dc}}}}$ respectively are the ac switching node current and voltage (i.e., before the output filter), $\bm{i}_{\tb\alpha\beta}$ and $\bm{v}_{\tb\alpha\beta}$ are the output current and voltage. 

To obtain a realistic model of the dc energy source, we model its response time by a first order system
\begin{equation}\label{eq:source_del}
    \tau_\text{\upshape{dc}}\dot{i}_\tau=i_\text{\upshape{dc}}^\star-i_\tau,    
\end{equation}
where $i^\star_{\text{\upshape{dc}}}$ is the dc current reference, $\tau_{\text{\upshape{dc}}}$ is the dc source time constant, and $i_\tau$ denotes the current provided by the dc source. Moreover, the dc source current limitation is modeled by the saturation function
\begin{equation}\label{eq:source_sat}
{i}_{\text{\upshape{dc}}}=\sat\left(i_\tau,i^{\text{dc}}_{\max}\right)=\begin{cases}
i_\tau&\If~|i_\tau|<|i^{\text{dc}}_{\max}|,\\    
\sgn\left(i_\tau\right)i^{\text{dc}}_{\max}&\If~|i_\tau|\geq|i^{\text{dc}}_{\max}|,
\end{cases}
\end{equation}
where $i^{\text{dc}}_{\max}$ is the maximum dc source current. Note that we implicitly assume that some storage element is present so that the dc source can support bidirectional power flow. We remark that the converter must limit its ac current to protect the semiconductor switches \cite{denis_migrate_2018}. Typically the current limitation is implemented via an underlying current controller (see Section \ref{GFCC}). 
}
{\tb
\subsection{Synchronous Machine Model}\label{syncmodel}
In this work we adopt an 8th order (i.e., including six electrical and two mechanical states), balanced, symmetrical, three-phase SM with a field winding and three damper windings on the rotor \cite[Fig. 3.1]{SP98}
\begin{subequations}\label{eqs:SM}
\begin{align}
\dot{\theta}&=\omega,\label{eqs:SM1}\\
J\dot{\omega}&=T_\text{m}-T_\text{e}-T_\text{f},\label{eqs:SM2}\\
\dot{\bm{\psi}}_\text{s,dq}&=\bm{v}_\text{s,dq}-r_\text{s}\bm{i}_\text{s,dq}-\mathcal{J}_2\bm{\psi}_\text{s,dq},\label{eqs:SM3}\\
\dot{\psi}_\text{f,d}&=v_\text{f,d}-r_\text{f,d}i_\text{f,d},\label{eqs:SM4}\\
\dot{\bm{\psi}}_\text{D}&=\bm{v}_\text{D}-\mathcal{R}_\text{D}\bm{i}_\text{D},\label{eqs:SM5}
\end{align}
\end{subequations}
where $\theta$ denotes the rotor angle, $J$ is the inertia constant, $\omega$ is the rotor speed, $T_\text{m}$, $T_\text{e}$ and $T_\text{fw}$ denote the mechanical torque, electrical torque, and the friction windage torque (see \cite[Sec. 5.7]{SP98}). Moreover, ${\bm{\psi}}_\text{s,dq}=[\psi_\text{s,d}~ \psi_\text{s,q}]^\top$, ${\bm{v}}_\text{s,dq}=[v_\text{s,d}~ v_\text{s,q}]^\top$, and ${\bm{i}}_\text{s,dq}=[i_\text{s,d}~ i_\text{s,q}]^\top$ denote the stator winding flux, voltage, and current in dq-coordinates (with angle $\theta$ as in \eqref{eqs:SM1}), $\mathcal{J}_2=[\begin{smallmatrix}
0&-1\\1&0
\end{smallmatrix}]$ denotes the $90^\circ$ rotation matrix, ${\psi}_\text{f,d}$, ${v}_\text{f,d}$, and ${i}_\text{f,d}$ denote the d-axis field winding flux, voltage and current. Furthermore, $r_\text{s}$ and $r_\text{f,d}$ denote the stator and field winding resistances, ${\bm{\psi}}_\text{D}=[{\psi}_\text{1d}~ {\psi}_\text{1q}~ {\psi}_\text{2q}]^\top, \bm{v}_\text{D}=[{v}_\text{1d}~ {v}_\text{1q}~ {v}_\text{2q}]^\top$ and $\bm{i}_\text{D}=[{i}_\text{1d}~ {i}_\text{1q}~ {i}_\text{2q}]^\top$ are the linkage flux, voltage and current associated with three damper windings and  $\mathcal{R}_\text{D}=\text{diag}\left({r}_\text{1d}, {r}_\text{1q}, {r}_\text{2q}\right)$ denotes the diagonal matrix of the damper winding resistances. Note that the friction term is conventionally expressed as a speed dependent term e.g., $T_\text{f}=D_\text{f}\omega$ \cite[Sec. 5.7]{SP98}. Furthermore, the damping torque associated with the damper windings is included in the SM model \eqref{eqs:SM} via the damper winding dynamics \eqref{eqs:SM5}. For more details on the SM modeling and parameters computation the reader is referred to \cite[Sec. 3.3]{SP98},\cite[Chap. 4]{kundur1994power}. 

We augment the system \eqref{eqs:SM} with a ST1A type excitation dynamics with built-in automatic voltage regulator (AVR) \cite[Fig. 21]{noauthor_ieee_2016}. To counteract the well-known destabilizing effect of the AVR on the synchronizing torque, we equip the system with a  simplified power system stabilizer comprised of a two-stage lead-lag compensator \cite[Sec. 12.5]{kundur1994power}.
Lastly, the governor and turbine dynamics are respectively modeled by proportional speed droop control and first order turbine dynamics
\begin{subequations}\label{Eq_pdel}
\begin{align}
p&=p^\star+{d_p}\left(\omega^\star-\omega\right), \label{Eq_sdroop}\\
\tau_{\text{g}}{\dot{p}}_{\tau}&={p}-{p_{\tau}},   
\end{align}
\end{subequations}
where $p^\star$ denotes the power set-point, $p$ is the governor output, $d_p$ denotes the governor speed droop gain, and $\omega^\star$ and $\omega$ denote nominal and measured frequency, respectively. Furthermore, $\tau_{\text{g}}$ denotes the turbine time constant and $p_\tau$ denotes the turbine output power. We refer the reader to \cite[Fig. 2]{Uros2019} for an illustration of the interplay between the SM model, the excitation dynamics, the PSS and governor dynamics. Lastly, the Matlab/Simulink implementation of the SM model can be found in \cite{model}.
}
\begin{figure}[b!]
    \centering
    \includegraphics[clip, trim=0.1cm 0cm 0cm 0cm,width=0.6\textwidth]{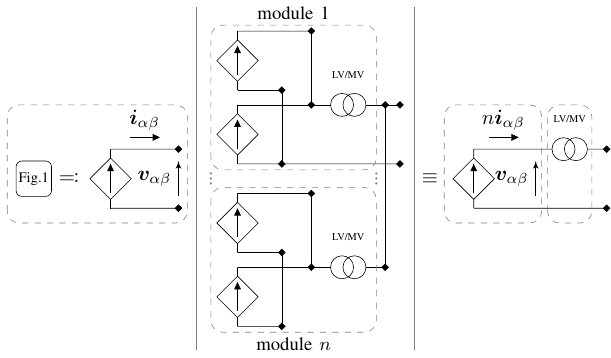}
    \caption{Equivalent model of an individual converter module (left), large-scale multi-converter system consisting of $n$ identical modules (middle), and aggregate model (right).\label{agg_model}}
\end{figure}
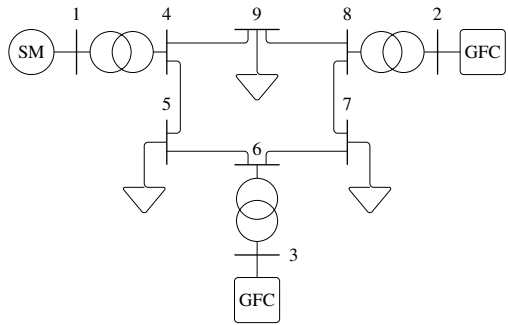
\begin{figure}[b!]
\centering
\tikzstyle{comp} = [coordinate]
\makeatletter
\begin{circuitikz}[scale=1.2,line width=0.1mm,rounded corners=0.5mm]
    {\fontfamily{ptm}\selectfont    
        \draw(0,1) circle(0.25cm);
        \draw(0.25,1)to[short]++(0.25,0);
        \draw(0.5,1)to[voosource]++(1,0);            
        \begin{scope}[line width=0.175mm]
            \draw (0.5,0.75)--node[pos=1,above]{\scriptsize 1}(0.5,1.25);
            \draw (1.5,0.75)--node[pos=1,above]{\scriptsize 4}(1.5,1.25);
            \draw (2.25,1.25)--node[pos=0.5,above]{\scriptsize 9}(2.75,1.25);
            \draw (1.5,-0.25)--node[pos=1,above]{\scriptsize 5}(1.5,0.25);
            \draw (2.25,-0.25)--node[pos=0.5,above]{\scriptsize 6}(2.75,-0.25);
            \draw (2.25,-1.25)--node[pos=1,right]{\scriptsize 3}(2.75,-1.25);
            \draw (3.5,0.75)--node[pos=1,above]{\scriptsize 8}(3.5,1.25);
            \draw (3.5,-0.25)--node[pos=1,above]{\scriptsize 7}(3.5,0.25);
            \draw (4.5,0.75)--node[pos=1,above]{\scriptsize 2}(4.5,1.25);
            \node[] at (0,1) {\scriptsize SM};
            \node[] at (5,1) {{\scriptsize GFC}};	
            \node[] at (2.5,-1.75) {\scriptsize GFC};    
        \end{scope}            
        \draw (1.5,1.1)-|(2.4,1.25);
        \draw (1.5,0.9)--(1.65,0.9)--(1.65,0.1)--(1.5,0.1);
        \draw (1.5,-0.1)-|(2.4,-0.25);
        \draw (2.6,-0.25)|-(3.5,-0.1);
        \draw (3.5,0.1)--(3.35,0.1)--(3.35,0.9)--(3.5,0.9);
        \draw (3.5,1.1)-|(2.6,1.25);
        \draw (2.5,1.25)--(2.5,0.75)--(2.75,0.75)--(2.5,0.45)--(2.25,0.75)--(2.5,0.75);
        \draw (1.5,0)--(1.25,0)--(1.25,-0.5)--(1.5,-0.5)--(1.25,-0.8)--(1,-0.5)--(1.25,-0.5);
        \draw (3.5,0)--(3.75,0)--(3.75,-0.5)--(4,-0.5)--(3.75,-0.8)--(3.5,-0.5)--(3.75,-0.5);            
        \draw(4.75,0.75) rectangle (5.25,1.25);
        \draw(4.75,1)to[short]++(-0.25,0);
        \draw(4.5,1)to[voosource]++(-1,0);
        \draw(2.5,-1.25)to[voosource]++(0,1);
        \draw(2.25,-2) rectangle (2.75,-1.5);
        \draw(2.5,-1.5)to[short]++(0,0.25);            
    }
\end{circuitikz}
    \caption{IEEE 9-bus test system with a synchronous machine, two large-scale multi-converter systems (i.e., aggregate GFCs), and constant impedance loads. \label{9bus_system}} 
\end{figure}
\begin{figure*}[t!]
    \centering
    \includegraphics[width=\textwidth]{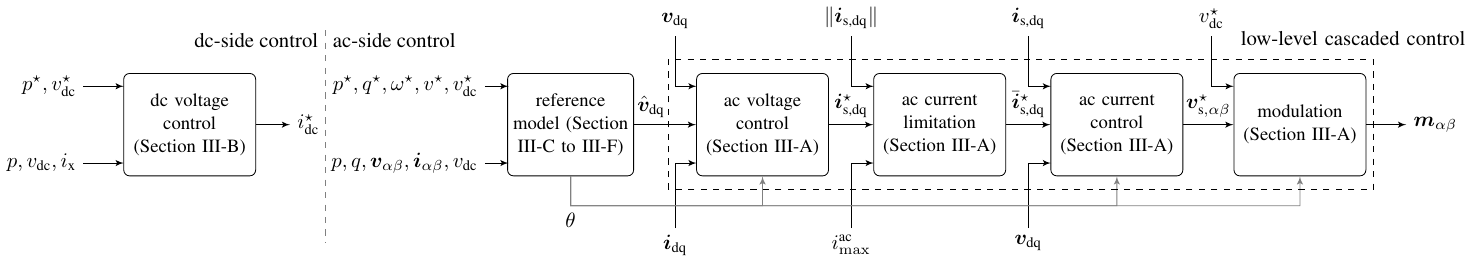}
    \caption{Grid-forming control architecture with reference models described in subsections \ref{subsec:droop} to \ref{subsec:dVOC}.\label{fig:architecture}}
\end{figure*}
\subsection{Network Model}
To study the transmission level dynamics of a low-inertia power system, we use Sim Power Systems to perform an EMT simulation of the IEEE 9-bus test system shown in Figure \ref{9bus_system} \cite{tayyebi_grid-forming_2018,ZMT11}. We model the lines via nominal $\pi$ sections (i.e., with RLC dynamics), model the transformers via three-phase linear transformer models, and consider constant impedance loads (see Table \ref{Table} for the parameters). We emphasize that the line dynamics cannot be neglected in the presence of grid-forming converters due to potential adverse interactions between their fast response and the line dynamics \cite{VHH+18,gros_effect_2018,Uros2019}.
\begin{remark}{(Aggregate Converter Model)}\\
\blue{\normalfont
In our case study, each GFC in Figure \ref{9bus_system} is an aggregate model of 200 commercial converter modules (see Table \ref{Table} for the parameters). Each module is rated at $500$ kVA and the aggregate model is rated at $100$ MVA, which is equal to the SM rating. Each module is interfaced to a medium voltage line via a LV-MV transformer (see Figure \ref{agg_model}). We derive the parameters of the aggregate transformer model by assuming a parallel connection of 100 commercial transformers rated at $1.6$ MVA (see Table \ref{Table}). A detailed presentation and derivation of the model aggregation is out of the scope of this work, but we follow developments analogous to \cite{purba_reduced-order_2017,Khan2018,PJRJBD19} in deriving the equivalent aggregate converter parameters.}
\end{remark} 

\section{Grid-Forming Control Architectures}\label{GFCC}
{Grid-forming control strategies control (see Figure \ref{fig:architecture}) a converter through the reference current $i^\star_\text{dc}$ for the dc energy source \eqref{eqs:model1} and the modulation signal $\bm{m}_{\tb\alpha\beta}$ for the dc-ac conversion stage \eqref{eqs:model} (see Figure \ref{fig:ConverterModel}). In the following, we briefly review the low-level cascaded control design (i.e., ac voltage control, current limitation and control) for two-level voltage source converters that tracks a voltage reference provided by a reference model (i.e., grid-forming control).} Moreover, we propose a controller for the converter dc voltage which defines the reference dc current. Because their design is independent of the choice of the reference model, we first discuss the cascaded voltage  / current control and the dc-side control. Subsequently, we review four common grid-forming control strategies. For each strategy, we describe the angle dynamics, frequency dynamics and ac voltage magnitude regulation. 
{\tb Throughout this section we will employ the three phase abc, $\alpha\beta$ and $\text{\upshape{dq}}$-coordinates (see \cite[Sec. 4.5 and 4.6]{YI10} for details on the transformations). We remind the reader that the Simulink implementation of the controls presented in the forthcoming subsections is available online \cite{model}.}
\subsection{Low-Level Cascade Control Design}\label{subsec:innercontr}
\subsubsection{AC Voltage Control}
we employ a standard converter control architecture that consists of a reference model providing a reference voltage $\hat{\bm{v}}_{\tb \text{dq}}$ with angle $\angle \hat{\bm{v}}_{\tb \text{dq}}= \theta$ and magnitude $\norm{\hat{\bm{v}}_{\tb \text{dq}}}$. The modulation signal $\bm{m}_{\tb\alpha\beta}$ is determined by cascaded  proportional-integral (PI) controllers that  are implemented  in $\text{\upshape{dq}}$-coordinates (rotating with  the reference angle $\theta$) and track the voltage reference $\hat{\bm{v}}_{\tb \text{dq}}$ (see \cite{PPG07}). The voltage tracking error $\hat{\bm{v}}_{\tb \text{dq}}-\bm{v}_{\tb \text{dq}}$ is used to compute the reference $\bm{i}^\star_{\tb\text{\upshape{s,dq}}} = [i^\star_{\text{\upshape{sd}}}~ i^\star_{{\text{\upshape{sq}}}}]^{\top}$ for the switching node current $\bm{i}_{\tb\text{\upshape{s,dq}}}$
\begin{subequations}\label{Eq_vloop}
\begin{align}
{\dot{\bm{x}}}_{v,{\tb\text{dq}}}&=\hat{\bm{v}}_{\tb \text{dq}}-\bm{v}_{\tb \text{dq}},
\\
\bm{i}^\star_{\tb\text{\upshape{s,dq}}} &\coloneqq \underbrace{\bm{i}_{\tb \text{dq}}+C\omega{\mathcal{J}_2}\bm{v}_{\tb \text{dq}}}_{\text{feed-forward terms}}+\underbrace{\mathcal{K}_{v,{\text{p}}}\left(\hat{\bm{v}}_{\tb \text{dq}}-\bm{v}_{\tb \text{dq}}\right)+\mathcal{K}_{v,{\text{i}}}{\bm{x}_{\tb v,\text{dq}}}}_{\text{PI control}}.\label{Eq_vloop:2}
\end{align}
\end{subequations}
Here {$\bm{x}_{v,{\tb\text{dq}}}=[x_{v,\text{\upshape{d}}}~x_{v,\text{\upshape{q}}}]^{\top}$} denotes the integrator state, ${\bm{v}}_{\tb \text{dq}}=[{v}_{\text{\upshape{d}}}~v_{\text{\upshape{q}}}]^{\top}$ denotes the output voltage measurement, $\hat{\bm{v}}_{\tb \text{dq}}=[\hat{v}_{\text{\upshape{d}}}~0]^{\top}$ denotes the reference voltage, $\bm{i}_{\tb \text{dq}}=[i_{\text{\upshape{d}}}~i_{\text{\upshape{q}}}]^{\top}$ denotes the output current, $\mathcal{I}_2$ is the 2-D identity matrix, $\mathcal{K}_{v,\text{p}}=k_{v,\text{p}}\mathcal{I}_2$ and $\mathcal{K}_{v,\text{i}}=k_{v,\text{i}}\mathcal{I}_2$ denote diagonal matrices of proportional and integral gains, respectively. 

\subsubsection{AC Current Limitation}
\blue{
We assume that the underlying current controller or low-level protections of the converter limit the ac current. We model this in abstraction by scaling down the reference current $\norm{\bm{i}^\star_{\text{\upshape{s,dq}}}}$ if it exceeds the pre-defined converter current limit $i^{\text{ac}}_{\max}$ \cite[Sec. III]{TWDB19}, i.e., 
\begin{align}\label{eq:ac_limiter}
\bar{\bm{i}}^\star_{\text{\upshape{s,dq}}}\coloneqq\begin{cases}
\bm{i}^\star_{\text{\upshape{s,dq}}}~&\If~\norm{\bm{i}_{\text{\upshape{s,dq}}}}\leq i^{\text{ac}}_{\max},\\    
\gamma_i\bm{i}^\star_{\tb\text{\upshape{s,dq}}}~&\If~\norm{\bm{i}_{\text{\upshape{s,dq}}}}> i^{\text{ac}}_{\max},
\end{cases}
\end{align}
where $\bar{\bm{i}}^\star_{\text{\upshape{s,dq}}}$ denotes the limited reference current that preserves the direction of $\bm{i}^\star_{\tb\text{\upshape{s,dq}}}$ and $\gamma_i\coloneqq\left({i^{\text{ac}}_{\max}}/{\norm{\bm{i}^\star_{\tb\text{\upshape{s,dq}}}}}\right)$. We emphasize that limiting the ac current can have a strong impact on the stability margins and dynamics of grid-forming power converters \cite{Linbin}. While numerous different ad-hoc strategies have been proposed to limit the ac current injection of voltage source converters with grid-forming controls  \cite{WLBL15,denis_migrate_2018,ZK17,GDS15,SHGM17,PD15,TWDB19} the problem of designing a robust ac current limitation strategy that effectively reacts to load-induced over-current and grid faults is an open research problem. Moreover, complex current limitation strategies typically require careful tuning of the controllers. Therefore, to provide a clear and concise comparison of existing grid-forming control solutions, we use the simple ac current limiting strategy \eqref{eq:ac_limiter}.}

\subsubsection{AC Current Control}
in order to implement this scheme, a PI controller for the current $\bm{i}_{\tb\text{\upshape{s,dq}}}=[i_{\text{\upshape{s,d}}}~i_{\text{\upshape{s,q}}}]^{\top}$ is used to track $\bar{\bm{i}}^\star_{\tb\text{\upshape{s,dq}}}$
\begin{subequations}\label{Eq_iloop}
\begin{align}
\dot{\bm{x}}_{\tb i,\text{\upshape{dq}}} &= \bar{\bm{i}}^\star_{\tb\text{\upshape{s,dq}}}-\bm{i}_{\tb\text{\upshape{s,dq}}},\\
\bm{v}^\star_{\tb\text{\upshape{s,dq}}}&\coloneqq\underbrace{\bm{v}_{\tb\text{\upshape{dq}}}+\left(L\omega {\mathcal{J}_2}+{R}\mathcal{I}_2\right)\bm{i}_{\tb\text{\upshape{s,dq}}}}_\text{feed-forward terms}+\underbrace{\mathcal{K}_{i,\text{p}}\left(\bm{i}^\star_{\tb\text{\upshape{s,dq}}}-\bm{i}_{\tb\text{\upshape{s,dq}}}\right)+\mathcal{K}_{i,\text{i}}{ \bm{x}}_{i,\text{\tb dq}}}_\text{PI control},\label{Eq_iloop:2}
\end{align}
\end{subequations}
where $\bm{v}^\star_{\tb\text{\upshape{s,dq}}} = [v^\star_{\text{\upshape{s,d}}}~ v^\star_{\text{\upshape{s,q}}}]^{\top}$ is the reference for the switching node voltage (i.e., before output filter in Figure \ref{fig:ConverterModel}), {$\bm{x}_{i,{\tb \text{dq}}}=[x_{i,\text{\upshape{d}}}~x_{i,\text{\upshape{q}}}]^{\top}$} denotes the integrator state, and ${\mathcal{K}}_{i,{\text{p}}}=k_{{\text{p},i}}\mathcal{I}_2$ and ${\mathcal{K}}_{i,{\text{i}}}=k_{i,{\text{i}}}\mathcal{I}_2$ denote the diagonal matrices of proportional and integral gains, respectively. Note that the first two terms of the right hand side of \eqref{Eq_vloop:2} and \eqref{Eq_iloop:2} are feed-forward terms. Finally, the modulation signal $\bm{m}_{\tb\alpha\beta}$ is given by 
\begin{align}\label{Eq_mod}
\bm{m}_{\tb\alpha\beta}=\frac{2\bm{v}^\star_{\tb\text{\upshape{s}},\alpha\beta}}{v^\star_{\text{\upshape{dc}}}},  
\end{align}
where $\bm{v}^\star_{\tb\text{\upshape{s}},\alpha\beta}$ is the $\alpha\beta$-coordinates image of $\bm{v}^\star_{\tb\text{\upshape{s,dq}}}$ defined in \eqref{Eq_iloop} and $v^\star_{{\text{\upshape{dc}}}}$ denotes the nominal converter dc voltage.
\begin{figure}[]
    \centering
    \includegraphics[clip, trim=0.85cm 0cm 0cm 0cm,width=0.7\textwidth]{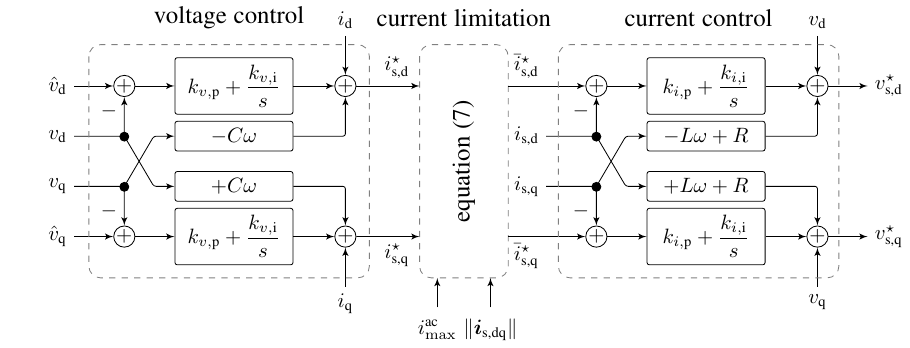}
    \centering    
    \caption{\tb Block diagram of the low-level cascaded control design \eqref{Eq_vloop}-\eqref{Eq_iloop} in $\text{\upshape{dq}}$-coordinates rotating with the angle $\theta$ provided by the reference model.\label{VC_loops}}
\end{figure}
\subsection{DC Voltage Control}\label{subsec:dccontr}
The dc current reference ${i}_{{\text{\upshape{dc}}}}^\star$ that is tracked by the controllable dc source \eqref{eq:source_del} is given by a proportional control for the dc voltage and feed-forward terms based on the nominal ac active power injection $p^\star$ and the filter losses
\begin{align}\label{Eq_idc_m}
{i}_{{\text{\upshape{dc}}}}^\star=\underbrace{k_{{\text{\upshape{dc}}}}\left(v_{{\text{\upshape{dc}}}}^\star-{v_{{\text{\upshape{dc}}}}}\right)}_{\text{proportional control}}+\underbrace{\frac{p^\star}{v_{{\text{\upshape{dc}}}}^\star}+\left(G_{{\text{\upshape{dc}}}}{v_{{\text{\upshape{dc}}}}}+\frac{v_{{\text{\upshape{dc}}}}i_{\text{\upshape{x}}}-{{p}}}{v_{{\text{\upshape{dc}}}}^\star}\right)}_{\text{power injection and loss feed-forward}},
\end{align}
where $v_{{\text{\upshape{dc}}}}i_{\text{\upshape{x}}}$ is the dc power flowing into the switches, $p$ is the ac power injected into the grid, and the last term on the right hand side of \eqref{Eq_idc_m} implements a feed-forward power control that compensates the filter losses. The loss compensation is required to ensure exact tracking of the power reference by matching control (see Section \ref{subsec:match}) and also improves the dc voltage regulation for the other control strategies considered in this study. Thus, to ensure a fair comparison, we apply \eqref{Eq_idc_m} for all control strategies discussed throughout this work.
\subsection{Droop Control}\label{subsec:droop}
Droop control resembles the speed droop property \eqref{Eq_sdroop} of the SM governor \cite{chandorkar_control_1993} and trades off deviations of the power injection (from its nominal value $p^\star$) and frequency deviations (from $\omega^\star$)
\begin{subequations}
\label{Eqtheta}
\begin{align}
{\dot{\theta}}&=\omega,\\
\omega&=\omega^\star+d_{\omega}\left(p^\star-p\right),\label{eq:droop_omega}
\end{align}
\end{subequations}
where $d_{\omega}$ denotes the droop gain. To replicate the service provided by the automatic voltage regulator (AVR) of synchronous machines we use a PI controller acting on the output voltage error
\begin{align}\label{Eq_vD}
    \hat{\bm{v}}_{\text{\upshape{d}}}=k_\text{p} \left(v^\star-\norm{\bm{v}_{\tb\text{dq}}}\right) + k_\text{i} \int_{0}^{t} \left({v}^\star-{\norm{\bm{v}_{\tb\text{dq}}(\tau)}}\right)\diff\tau.
\end{align}
to obtain the direct axis reference $\hat{v}_{\text{\upshape{d}}}$ for the underlying voltage loop ($v^\star$ and ${\norm{\bm{v}_{\tb\text{dq}}}}$ are the reference and measured voltage magnitude). We remark that {$\hat{v}_\text{q}=0$ and} the reactive power injection varies such that exact voltage regulation is achieved. 
\begin{figure}[b!!]
    \centering
    \includegraphics[clip, trim=0.65cm 0cm 0cm 0cm,width=0.49\textwidth]{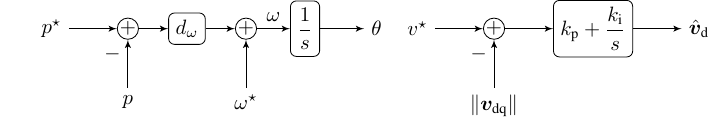}
    \hfill
    \includegraphics[clip, trim=0.65cm 0cm 0cm 0cm,width=0.49\textwidth]{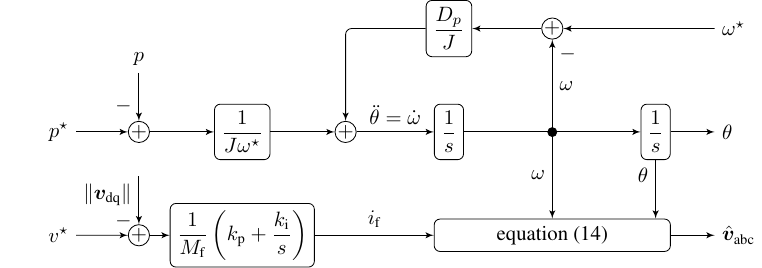}
    \caption{Droop control frequency and ac voltage control block diagrams based on \eqref{Eqtheta} and \eqref{Eq_vD} (left), block diagram of a grid-forming VSM based \eqref{Eq_Theta_Omega_Sync}-\eqref{EQ_IF} (right).}\label{DroopControl}
\end{figure}
\subsection{Virtual Synchronous Machine}\label{subsec:VSM}
Many variations of virtual synchronous machines (VSMs) have been proposed \cite{darco_virtual_2013,zhong_synchronverters:_2011}. In this work, we consider the frequency dynamics induced by the synchronverter \cite{zhong_synchronverters:_2011} 
\begin{subequations}\label{Eq_Theta_Omega_Sync}
\begin{align}
{\dot{\theta}}&=\omega,\\
{J}_\text{r}\dot{\omega}&=\frac{1}{\omega^\star}\left(p^\star-p\right)+{D_{p}}\left(\omega^\star-\omega\right)\label{eq:omega_vsm},
\end{align}
\end{subequations}
{\tb where $D_{p}$ is commonly referred to as damping factor \cite[Sec. II-b]{zhong_synchronverters:_2011} that, in abstraction, models the governor speed droop coefficient}, $J_\text{r}$ is the virtual rotor's inertia constant. Note that the dynamics \eqref{Eq_Theta_Omega_Sync} reduce to droop control \eqref{Eqtheta} when using $J_\text{r}/D_{p}\approx0$ as recommended in \cite{zhong_synchronverters:_2011}. These angle dynamics capture the main salient features of virtual synchronous machines, but do not suffer from drawbacks of more complicated implementations (see \cite{tayyebi_grid-forming_2018} for a discussion). The three-phase voltage induced by the VSM is given by
\begin{align}
\hat{\bm{v}}_{\tb\text{\upshape{abc}}}=2\omega M_\text{f}i_\text{f}\begin{bmatrix}\sin {(\theta)} & \sin{\left(\theta-{\tfrac{2\pi}{3}}\right)} & \sin{\left(\theta-{\tfrac{4\pi}{3}}\right)}\end{bmatrix}^{\top}
\label{EQ_VS},
\end{align}
where $M_\text{f}$ and $i_\text{f}$ are respectively the virtual mutual inductance magnitude and excitation current. Similar to \eqref{Eq_vD}, we utilize input $i_\text{f}$ to achieve exact ac voltage regulation via PI control and thereby replicate the function of the AVR of a synchronous machine
\begin{align}\label{EQ_IF}
i_{\text{f}}=\frac{k_\text{p}}{M_\text{f}} \left({v}^\star-{\norm{\bm{v}_{\tb\text{dq}}}}\right) + \frac{k_\text{i}}{M_\text{f}} \int_{0}^{t} \left({v}^\star-{ \norm{\bm{v}_{\tb\text{dq}}(\tau)}}\right) \diff \tau.
\end{align}
Transforming $\hat{\bm{v}}_{\tb\text{abc}}$ to $\text{\upshape{dq}}$-coordinates with $\theta$ and $\omega$ as in \eqref{Eq_Theta_Omega_Sync}, the voltage and current loops and modulation signal generation remain the same as \eqref{Eq_vloop}--\eqref{Eq_mod}.
\subsection{Matching Control}\label{subsec:match}
\begin{figure}[t!]
        \centering
        \includegraphics[clip, trim=0.65cm 0cm 0cm 0cm,width=0.4\textwidth]{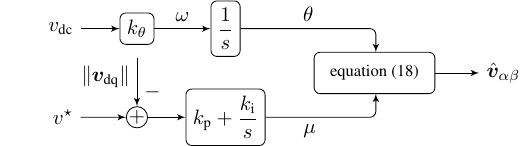}    
    \hfill
        \centering
        \caption{Matching control block diagram based on \eqref{Eq_theta_match}-\eqref{EQ_modM}.}\label{Matching}
\end{figure}
Matching control is a grid-forming control strategy that exploits structural similarities between power converters and SMs \cite{cvetkovic_modeling_2015,arghir_grid-forming_2017,curi_control_2017,Catalin2019,Linbin2017} and is based on the observation that the dc-link voltage - similar to the synchronous machine frequency - indicates power imbalances. Hence, the dc voltage, up to a constant factor, is used to drive the converter frequency. This control technique structurally matches the differential equations of a converter to those of a SM. Furthermore, analogous to the machine input torque, the dc current is used to control the ac power. The angle dynamics of matching control are represented by 
\begin{align}\label{Eq_theta_match}
    \dot{\theta}&= k_\theta v_{{\text{\upshape{dc}}}},
\end{align}
where $k_\theta\coloneqq\omega^\star/v_{{\text{\upshape{dc}}}}^\star$. Finally, the ac voltage magnitude is controlled through the modulation magnitude $\mu$ by a PI controller    
\begin{align}
\mu=k_\text{p} \left({v}^\star-{\norm{\bm{v}_{\tb\text{dq}}}}\right) + k_\text{i} \int_{0}^{t} \left({v}^\star-{\norm{\bm{v}_{\tb\text{dq}}(\tau)}}\right)\diff\tau.
\label{EQ_VM}
\end{align}
The reference voltage for the voltage controller in $\alpha\beta$-coordinates is given by:
    \begin{align}
    \hat{\bm{v}}_{\tb{\alpha\beta}}=\mu[-\sin{\theta}~\cos{\theta}]^{\top}.
    \label{EQ_modM}
    \end{align}
Transforming $\hat{\bm{v}}_{\tb{\alpha\beta}}$ to $\text{\upshape{dq}}$-coordinates with $\theta$ and $\omega$ as in \eqref{Eq_theta_match}, the voltage and current loops and modulation signal generation remain the same as \eqref{Eq_vloop}--\eqref{Eq_mod}. 

{\tb To further explain the matching concept, we replace $v_\text{dc}$ in \eqref{eqs:model1} by $\omega/k_\theta$ from \eqref{Eq_theta_match} resulting in
    \begin{subequations}\label{eqs:MatchSwing}    
        \begin{align}
        \dot{\theta}&=\omega,\\
        C_{{\text{\upshape{dc}}}}{\dot{\omega}}&=k_\theta{i}_{{\text{\upshape{dc}}}}-k_\theta i_{\text{\upshape{x}}}-G_{{\text{\upshape{dc}}}}\omega.\label{eq:vdc_matching}
        \end{align} 
    \end{subequations}
    Recalling the SM's angle and frequency dynamics \eqref{eqs:SM1}-\eqref{eqs:SM2} and replacing $T_\text{f}$ by $D_\text{f}\omega$ 
    \begin{subequations}\label{eqs:SMSwing}    
        \begin{align}
        \dot{\theta}&=\omega,\\
        J{\dot{\omega}}&=T_\text{m}-T_\text{e}-D_\text{f}\omega, \label{eqs:SMSwing:freq}
        \end{align} 
    \end{subequations} 
Comparing \eqref{eqs:MatchSwing} and \eqref{eqs:SMSwing} reveals the structural matching of GFC dynamics to that of SM.  Dividing \eqref{eq:vdc_matching} by $k^2_{\theta}$ to obtain the same units as in \eqref{eqs:SMSwing:freq} and matching variables results in $J_\text{r}=C_{{\text{\upshape{dc}}}}/k^2_{\theta}$, $D_\text{f}=G_{{\text{\upshape{dc}}}}/k^2_{\theta}$, $T_\text{m}=i_{{\text{\upshape{dc}}}}/k_{\theta}$, and $T_\text{e}=i_{\text{\upshape{x}}}/k_{\theta}$. In other words, using matching control the inertia constant of the GFC is linked to its internal energy storage, the dc-side losses $G_\text{dc}\omega$ are linked to the machine friction losses $D_\text{f}\omega$, and the frequency droop response is provided through the proportional dc voltage control $k_\text{dc}\left(v^\star_\text{dc}-v_\text{dc}\right)=\left(k_\text{dc}/k_\theta\right)\left(\omega^\star-\omega\right)$ (cf. \eqref{Eq_idc_m}). The structural matching induced by \eqref{Eq_theta_match} also extends to the converter ac filter and generator stator dynamics (see \cite{arghir_grid-forming_2017,Catalin2019} for a detailed derivation).
}
\subsection{Dispatchable Virtual Oscillator Control}\label{subsec:dVOC}
\begin{figure}[b!!]
        \centering   
        \includegraphics[clip, trim=0.65cm 0cm 0cm 0cm,width=0.49\textwidth]{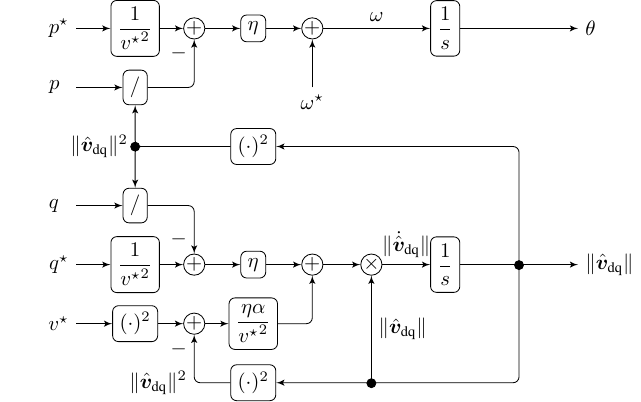}
        \centering
        \caption{Block diagram of dVOC in polar coordinates corresponding to \eqref{eq.dvoc.droop}. Note that singularity at $\norm{\hat{\bm{v}}_{\tb\text{dq}}}=0$ only appears in the dVOC implementation in polar coordinates but not in the implementation in rectangular coordinates, i.e., \eqref{Eq_dvoc_v}.\label{dVOC}}
\end{figure}
Dispatchable virtual oscillator control (dVOC) \cite{colombino_global_2017,gros_effect_2018,SG-SI-BJ-MC-DG-FD:18} is a decentralized grid-forming control strategy that guarantees almost global asymptotic stability for interconnected GFCs with respect to nominal voltage and power set-points \cite{colombino_global_2017,gros_effect_2018}. The analytic stability conditions for dVOC explicitly quantify the achievable performance and include the dynamics and transfer capacity of the transmission network \cite{gros_effect_2018}.  

The dynamics of dVOC in $\alpha\beta$-coordinates are represented by
\begin{align}\label{Eq_dvoc_v}
\!\!\dot{\hat{\bm{v}}}_{\tb\alpha\beta}=\omega^\star \mathcal{J}_2 \hat{\bm{v}}_{\tb\alpha\beta}+\eta\left( \mathcal{K}\hat{\bm{v}}_{\tb\alpha\beta}-\mathcal{R}_2(\kappa)\bm{i}_{\tb\alpha\beta}+ \frac{\alpha}{v^{\star2}} \left(v^{\star2}-\norm{\hat{\bm{v}}^2_{\tb\alpha\beta}}\right) \hat{\bm{v}}_{\tb\alpha\beta}\right)\!,
\end{align}
where $\hat{\bm{v}}_{\tb\alpha\beta}=[\hat{v}_\alpha~\hat{v}_\beta]^{\top}$ is the reference voltage, $\bm{i}_{\tb\alpha\beta}=[i_\alpha~i_\beta]^{\top}$ is current injection of the converter, the angle $\kappa\coloneqq\tan^{-1}\left( l\omega^\star/r\right)$ models the network inductance to resistance ratio, and $\eta$, $\alpha$ are positive control gains. Furthermore we have
\begin{align*}
\mathcal{R}_2(\kappa)\coloneqq\begin{bmatrix}\cos{\kappa} & -\sin{\kappa}\\\sin{\kappa} & \cos{\kappa}\end{bmatrix},\mathcal{K}\coloneqq\frac{1}{{{v}^\star}^2}\mathcal{R}_2(\kappa)\begin{bmatrix}p^\star & q^\star\\-q^\star & p^\star\end{bmatrix},    
\end{align*}
where $\mathcal{R}_2(\kappa)$ is the 2-D rotation by $\kappa$. As shown in \cite{gros_effect_2018} the dynamics \eqref{Eq_dvoc_v} reduce to a harmonic oscillator if phase synchronization is achieved (i.e., $\mathcal{K}\hat{\bm{v}}_{\tb\alpha\beta}-\mathcal{R}_2(\kappa)\bm{i}_{\tb\alpha\beta}=0$) and $\norm{\hat{\bm{v}}_{\tb\alpha\beta}} = v^\star$ (i.e., $\left(v^{\star2}-\norm{\hat{\bm{v}}_{\tb\alpha\beta}}^2\right) \hat{\bm{v}}_{\tb\alpha\beta}=0$). Rewriting \eqref{Eq_dvoc_v} in polar coordinates for an inductive network (i.e., $\kappa = \pi/2$) reveals the droop characteristics (see \cite{colombino_global_2017,gros_effect_2018,SG-SI-BJ-MC-DG-FD:18}) of dVOC as 
\begin{subequations}\label{eq.dvoc.droop}
    \begin{align}
    \!\!\! \dot{\theta}&= \omega^\star + \eta \left(\frac{p^\star}{v^{\star2}} -\frac{p}{\norm{\hat{\bm{v}}_{\tb\text{dq}}}^2}\right),\label{eq.droop.p.approx} \\ 			
    \!\!\! \norm{\text{\mbox{$\dot{\hat{\bm{v}}}$}}_{\tb\text{dq}}}\;  & =\eta \left(\frac{q^\star}{{v}^{\star2}} - \frac{q}{\norm{\hat{\bm{v}}_{\tb\text{dq}}}^2} \right)\norm{\hat{\bm{v}}_{\tb\text{dq}}}+\frac{\eta \alpha}{v^{\star2}} \left( v^{\star2} - \norm{\hat{\bm{v}}_{\tb\text{dq}}}^2 \right) \norm{\hat{\bm{v}}_{\tb\text{dq}}}. \label{eq.droop.q.approx}			
    \end{align}
\end{subequations}
In other words, for a high voltage network and near the nominal steady state (i.e., $\norm{\hat{\bm{v}}_{\tb\text{dq}}} \approx v^{\star}$) the relationship between frequency and active power resemble that of standard droop control given in \eqref{Eqtheta} with $d_\omega = \eta / v^{\star2}$. Moreover, when choosing the control gain $\alpha$ to obtain post-fault voltages consistent with the other control algorithms described above, the first term in \eqref{eq.droop.q.approx} is negligible and \eqref{eq.droop.q.approx} reduces to the voltage regulator $\norm{\text{\mbox{$\dot{\hat{\bm{v}}}$}}_{\tb\text{dq}}}\; \approx - 2 \eta \alpha \left(\norm{\hat{\bm{v}}_{\tb\text{dq}}} - v^\star\right)$ near the nominal steady state.
\section{Network Case Study}\label{simulation}
In this section we present a comparison study of grid-forming control techniques in the presence of synchronous machines. In the forthcoming discussion, we use the test system shown in Figure \ref{9bus_system}. The parameters and control gains are given in Table \ref{Table}. The implementation in Simulink is available online \cite{model}. 

In order to avoid the delay associated with the frequency measurement and signal processing  introduced by standard synchronous reference frame phase-locked loop (SRF-PLL)\cite{poolla_placement_2018}, we use the mechanical frequency of the SM at node 1 in Figure \ref{9bus_system} to evaluate the post-disturbance system frequency (e.g., in Figures \ref{all_rcf}-\ref{all_mfd}). For the grid-forming converters we use the internal controller frequencies defined by \eqref{Eqtheta}, \eqref{Eq_Theta_Omega_Sync}, \eqref{Eq_theta_match} and \eqref{eq.dvoc.droop}. We remark that, in a real-world system and in an EMT simulation (in contrast to a phasor simulation), there is no well-defined frequency at the voltage nodes during transients, whereas the internal frequencies of the grid-forming units are always well-defined \cite[Sec. II-J]{milano_foundations_2018},\cite{Milano2017}. Lastly, we note that in all the forthcoming case studies we assume that the system is in steady-state at $t=0$.
\begin{figure}[b!]
   \centering
   \includegraphics[width=0.45\textwidth]{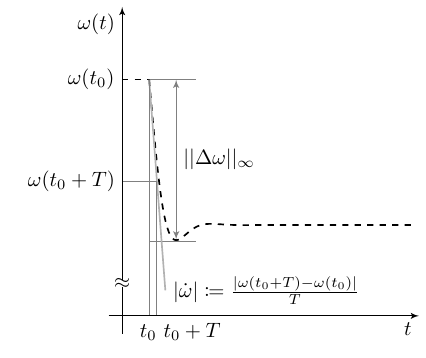}
   \caption{Post-event frequency nadir and RoCoF.}\label{fig:metrics}
\end{figure}
\subsection{Performance Metrics}
We adopt the standard power system frequency performance metrics i.e., maximum frequency deviation $||\Delta\omega||_\infty$ (i.e., frequency nadir/zenith) and RoCoF $|{\dot{\omega}}|$ (i.e., the slope of line tangent to the post-event frequency trajectory) defined by
\begin{subequations}\label{fmetrics}
    \begin{align}
    ||\Delta\omega||_\infty &\coloneqq \max_{t\geq t_0} |\omega^\star-\omega(t)|,\\
    |{\dot{\omega}}| &\coloneqq\frac{|\omega(t_0+T)-\omega(t_0)|}{T},
    \end{align}   
\end{subequations}
where $t_0>0$ is the time when the disturbance is applied to the system, and $T>0$ is the RoCoF calculation window \cite{milano_foundations_2018,poolla_placement_2018}. See Figure \ref{fig:metrics} for visual representation of the metrics described by \eqref{fmetrics}.
 In this work, we use $T=250 \mathrm{ms}$, which is in line with values suggested for protection schemes (see \cite[Table 1]{entsoefreq_2018}). Dividing the metrics \eqref{fmetrics} by the size of the magnitude of the disturbance results in a measure of the system disturbance amplification.
\subsection{Test Network Configuration and Tuning Criteria}\label{subsec:tuning}
In order to study the performance of the control approaches introduced in Section \ref{GFCC}, we apply the same strategy (with identical tuning) for both converters (i.e., at nodes 2 and 3 in Figure \ref{9bus_system}), resulting in four different SM-GFC paired models. As a benchmark, we also consider an all-SMs system with three identical SMs (i.e., at nodes 1-3). Selecting fair tuning criteria for the different control strategies is a challenging task. For this study, we tune the control parameters such that all generation units exhibit identical proportional load sharing behavior. Appendix \ref{app:tunning} presents our tuning criteria and derivation of some control parameters. 
\subsection{Impact of Grid-Forming Control on Frequency Metrics}\label{simulation_freq}
\begin{figure*}[b!!!]
    \centering
    \includegraphics[clip, trim=0.9cm 0.55cm 0.55cm 0.55cm,width=0.94\textwidth]{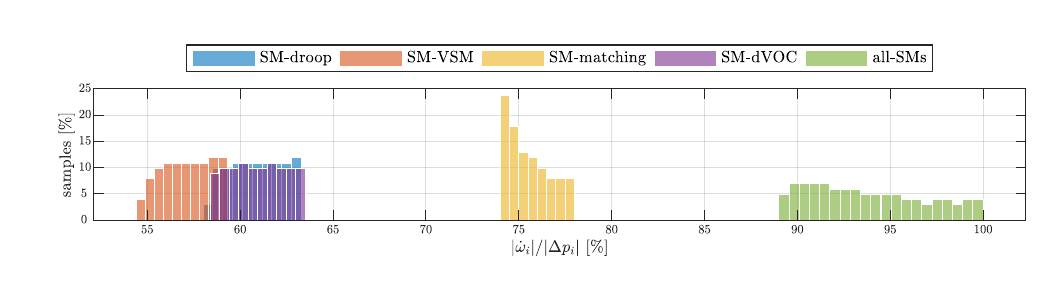}
    \vspace{-2mm}
    \caption{Normalized distribution of the RoCoF $|{\dot{\omega}}_i| / |\Delta p_i|$ of the synchronous machine frequency at node 1 for load disturbances $\Delta p_i$ ranging from $0.2$ p.u. to $0.9$ p.u. at node 7. For each load disturbance, $|{\dot{\omega}}_i| / |\Delta p_i|$ is normalized by the maximum value corresponding to the all-SMs configuration.\label{all_rcf}} 
    \centering
    \includegraphics[clip, trim=0.55cm 0.55cm 0.55cm 0.0cm,width=0.75\textwidth]{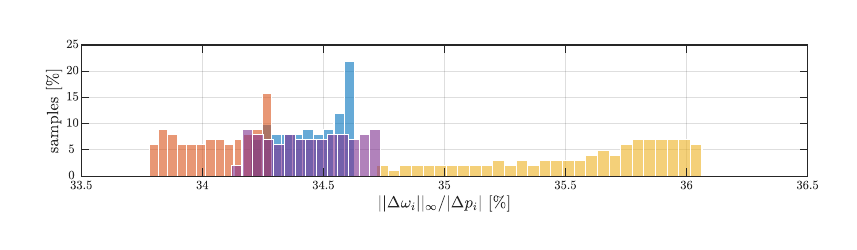}
    \includegraphics[clip,trim=0.55cm 0.55cm 0.55cm 0.0cm,width=0.19\textwidth]{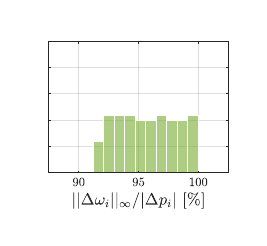}
    \vspace{-2mm}
    \caption{Normalized distribution of the nadir $||\Delta\omega_i||_{\infty} / |\Delta p_i|$ of the synchronous machine frequency at node 1 for load disturbances $\Delta p_i$ ranging from $0.2$ p.u. to $0.9$ p.u. at node 7. For each load disturbance, $||\Delta\omega_i||_{\infty} / |\Delta p_i|$ is normalized by the maximum value corresponding to the all-SMs configuration.    \label{all_mfd}}    
\end{figure*}

In this section we test the system behavior for different load disturbances $\Delta p_i$. The network base load $p_l$ is constant and uniformly distributed between nodes 5, 7 and 9 while $\Delta p_i$ is only applied at node 7. For each disturbance input we calculate $||\Delta\omega_i||_\infty$ and $|{\dot{\omega}}_i|$ for the SM at node 1 and normalize these quantities by dividing by $|\Delta p_i|$. Figures \ref{all_rcf} and \ref{all_mfd} illustrate the distribution of system disturbance input/output gains associated with introduced frequency performance metrics. Note that the network base load $p_l$ is 2 pu and the elements of the load disturbance sequence $\Delta p_i\in[0.2,0.9],~i=1,\dots,100$ are uniformly increasing by 0.007 pu starting from $p_1=0.2$ pu.

\begin{figure}[b!!]
        \centering
        \includegraphics[clip, trim=0.65cm 0.65cm 0.65cm 0.65cm,width=0.49\textwidth]{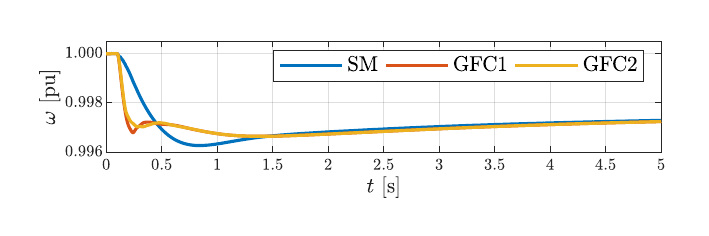}    
        \caption{Frequency of the system with two VSMs after a 0.75 pu load increase. The converters quickly synchronize with each other and then slowly synchronize with the machine.    \label{sync_f_stable}}
\end{figure}

Figures \ref{all_rcf} and \ref{all_mfd} suggest that, regardless of the choice of control strategy, the presence of grid-forming converters improves the metrics compared to the all-SM system. This possibly observation can be explained by the fast response of converters compared to the slow turbine dynamics, i.e., $\tau_\text{g}$ in \eqref{Eq_pdel} is larger than $\tau_{{\text{\upshape{dc}}}}$ in \eqref{eq:source_del}. Because of this, the converters reach frequency synchronization at a faster time-scale and then synchronize with the SM (see Figure \ref{sync_f_stable}). Overall, for any given disturbance input, the converters are able to react faster than the SM and the remaining power imbalance affecting the SM is smaller than in the all-SM system. This result highlights that the fast response of GFCs should be exploited instead of designing the controls of a converter (fast physical system) to emulate the slow response of synchronous machines \cite{milano_foundations_2018}.

We observe that droop control and dVOC exhibit very similar performance confirming the droop-like behavior of dVOC in predominately inductive networks (see \eqref{eq.dvoc.droop}). Moreover, the difference between droop control and VSM arises from the inertial (derivative control) term in \eqref{Eq_Theta_Omega_Sync} and the RoCoF is considerably higher when using matching control. This can be explained by the fact that VSM, droop control, and dVOC ignore the dc voltage and aggressively regulate ac quantities to reach angle synchronization, thus requiring higher transient peaks in dc current to stabilize the dc voltage (see Section \ref{subsec:instability}). Although improving RoCoF, this approach can lead to instability if the converter is working close to the rated power of the dc source, as shown in the next section. On the other hand, matching control regulates the dc link voltage both through the dc source and by adjusting its ac signals.

we selected the RoCoF calculation window according to the guideline \cite{entsoefreq_2018}, which accounts for noise and possible oscillations in the frequency signal. However, these guidelines were derived for a power system fully operated based on synchronous machines. Given that grid-forming converters introduce faster dynamics, machines are expected to reach the frequency nadir faster. Hence, a smaller RoCoF windows might need to be considered in a low-inertia power system to properly assess system performance. We note that the performance of the different grid-forming control strategies shown in Figure \ref{all_rcf} and \ref{all_mfd} is sensitive to the tuning of control gains and choice of RoCoF computation window. However, due to the comparably slow response of conventional generation technology the performance improvements for the system with grid-forming converters over the all-SM system persist for a wide range of parameters. Moreover, using comparable tuning (see Section \ref{subsec:tuning}) the differences between the different grid-forming techniques observed in this section are expected to remain the same.

\subsection{Instability Behavior -- Large Load Disturbance}\label{subsec:instability}
In this subsection we analyze the response of the grid-forming converters to large disturbances when the dc source is working close to its maximum rated values. In this case study, the dc-side current limitation of GFCs has a major impact on the overall system behavior. {\tb For the clarity of representation we first focus exclusively on the influence of dc source current limitation and neglect the ac current limitation. A test case highlighting the consequences of combined dc and ac limitation is presented in the next subsection (see the previous subsection).} 

\begin{figure}[]
\centering
\includegraphics[clip, trim=0.65cm 0.65cm 0.65cm 0.55cm,width=0.49\textwidth]{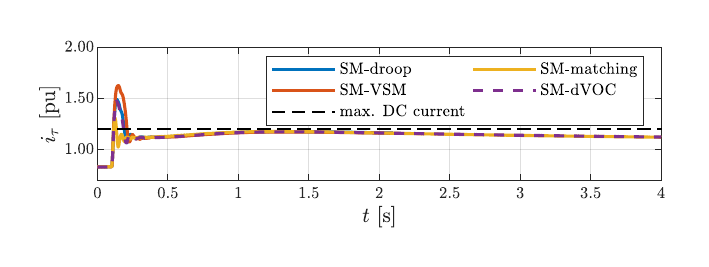}\hfill
\includegraphics[clip, trim=0.65cm 0.65cm 0.65cm 0.55cm,width=0.49\textwidth]{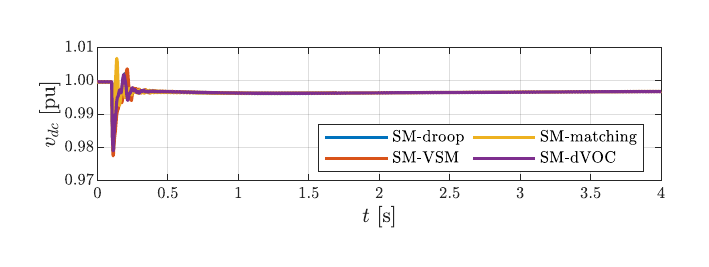}
\caption{dc current demand of the converter at node 2 ({\tb left}) and its dc voltage ({\tb right}) after a $0.75$ pu load disturbance.\label{fig:stable_idc}}
\centering
\includegraphics[clip, trim=0.65cm 0.65cm 0.65cm 0cm,width=0.49\textwidth]{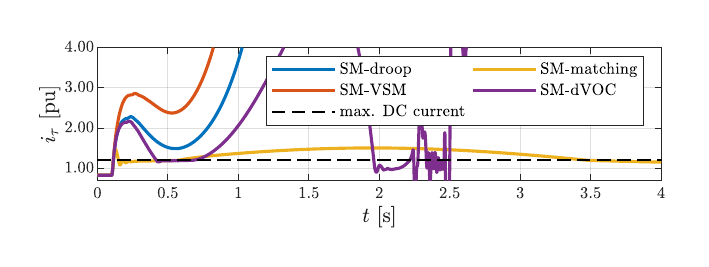}\hfill
\includegraphics[clip, trim=0.65cm 0.65cm 0.65cm 0cm,width=0.49\textwidth]{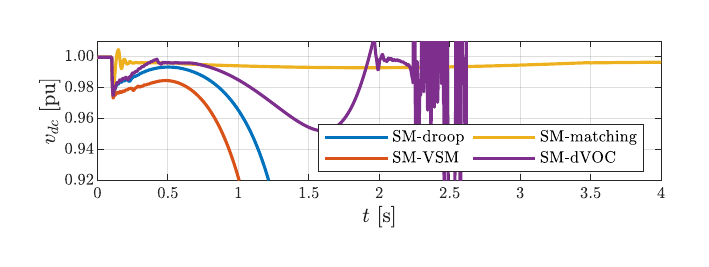}
\caption{dc current demand of the converter at node 2 ({\tb left}) and its dc voltage ({\tb right}) after a $0.9$ pu load disturbance.\label{fig:unstable_idc}}
\centering
\includegraphics[clip, trim=0.65cm 0.65cm 0.65cm 0cm,width=0.49\textwidth]{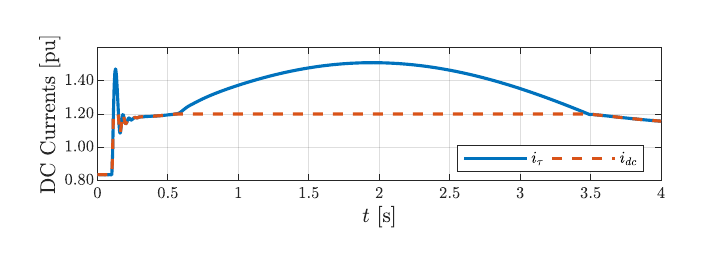}\hfill
\includegraphics[clip, trim=0.65cm 0.65cm 0.65cm 0cm,width=0.49\textwidth]{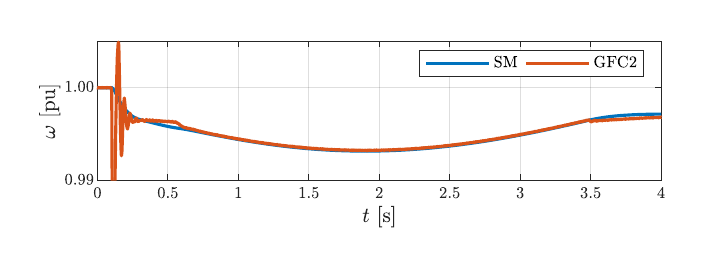}
\caption{dc current demand and saturated dc current (left), frequency of the converter (using matching control) at node 2 and SM after a $0.9$ pu load disturbance (right).\label{fig:frequency_matching}}
\includegraphics[clip, trim=0.65cm 0.65cm 0.65cm 0cm,width=0.49\textwidth]{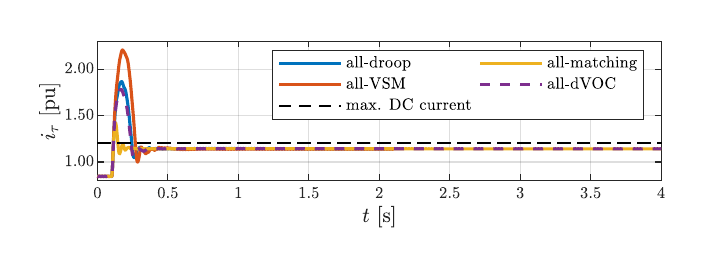}\hfill
\includegraphics[clip, trim=0.65cm 0.65cm 0.65cm 0cm,width=0.49\textwidth]{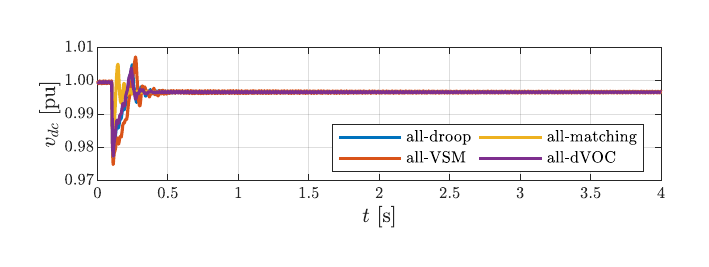}
\caption{dc current demand ({\tb left}) and dc voltage ({\tb right}) after a $0.9$ pu load disturbance in an all-GFC system.\label{fig:stable_idc_all_gfc}}
\end{figure}

To begin with, we set the network base load $p_l$ and load change {\tb $\Delta p$} to $2.25$ and $0.75$ pu respectively. Figure \ref{fig:stable_idc} shows the dc voltage and delayed dc current (i.e., before the saturation \eqref{eq:source_sat}) of the converter at node 2. Considering the base loading and large disturbance magnitude (which rarely occurs in transmission systems), it is interesting to observe that all controls remain stable. 

We now increase $\Delta p$ to $0.9$ pu (i.e., a total network load after the disturbance of 3.15 pu) which is equally shared by the SM and the GFCs. We expect a post-disturbance converter power injection of $1.05$ pu and $i_{{\text{\upshape{dc}}}}$ close to the dc current limit $i^{\text{dc}}_{\max}=\pm 1.2$ pu. Figure \ref{fig:unstable_idc} shows the dc voltage and delayed dc current before saturation for the converter at node 2 for a load increase of $\Delta p=0.9$ pu. For sufficiently large disturbance magnitude, $i_{{\text{\upshape{dc}}}}^\star$ and consequently $i_{\tau}$ exceed the current limit, i.e., the dc current $i_{{\text{\upshape{dc}}}}=i^{\text{dc}}_{\max}$ is saturated. We observe that the VSM, droop control, and dVOC fail to converge to a stable post-event equilibrium if the dc source is saturated for a prolonged time. In contrast, the GFCs controlled by matching reach a stable post-event equilibrium.

Pinpointing the underlying cause of instability is not straightforward. For the VSM and converters controlled by droop control and dVOC the dc-link capacitor discharges to provide $i_{\tau}-i^{\text{dc}}_{\max}$ (i.e., the portion of current demand which is not provided due to the saturation \eqref{eq:source_sat}). For droop control and the VSM this results in a collapse of the dc voltage. When the dc voltage drops below a certain limit - tripping the converter in practice - the ac voltage collapses and consequently the reactive power diverges far off the practical limits. This is followed by a similar instability behavior for active power. Because dVOC decreases its power injection if the voltage is low (i.e., $\norm{\bm{v}_\text{\tb dq}}^2 < v^{\star2}$ and $p$ and $q$ resulting from \eqref{Eq_dvoc_v} are lower) the dc voltage does not collapse immediately. Nonetheless, dVOC cannot reach a synchronous solution  and becomes unstable because the dc current is limited. 

Reducing the ac voltage control gains in \eqref{Eq_vD} and \eqref{EQ_VS}, i.e., enforcing stronger dc-ac dynamics time-scale separation, stabilizes all GFCs but less accurate tracking of the reference voltage $\hat{\bm{v}}_\text{\tb dq}$. Additionally, it has been observed that removing the inner voltage and current tracking loops \eqref{Eq_vloop}-\eqref{Eq_iloop}, i.e., using $\hat{\bm{v}}_\text{\tb dq}$ directly to obtain the modulation, stabilizes all three strategies. This supports the argument made in \cite{synchronverter2014} that removing the inner loops increases the bandwidth of the controller and can result in increased performance. However, we remark that these ad hoc remedies by no means guarantee stability for larger disturbances or different converter and network parameters. Finally, it has been reported that limiting the ac current also destabilizes droop control \cite{Linbin}.

{\tb In the next subsection we show that the instability GFCs controlled by control droop, VSMs and dVOC can be resolved by a suitable ac current limitation strategy that prevents depleting the dc voltage.}  
In contrast to the other techniques matching control succeeds to stabilize dc voltage despite saturation of the dc source. From a circuit-theoretic point of view this is only possible if the sum of the ac power injection and filter losses equals the approximately constant dc power inflow $v_{{\text{\upshape{dc}}}} i^{\text{dc}}_{\max}$. The converter can only inject constant ac power into the network if its angle difference with respect to the remaining devices in the network is constant. In the presence of the slow SM angle and frequency dynamics this implies that the GFCs need to synchronize to the SM (and with each other) so that the relative angle ${\theta}_{\text{GFC}}-\theta_{\text{SM}}=\theta_{\max}$ is constant.  

This synchronization is achieved through the dc voltage imbalance, i.e., as long as the dc voltage deviates from its nominal value matching control adjusts its voltage angle (see \eqref{Eq_theta_match}). In particular, the brief initial frequency transient (after the dc current reaches its limit) shown in Figure \ref{fig:frequency_matching} balances the power flowing in and out of the dc capacitor and results in an angle difference to the SM of $\theta_{\max}$. Overall, this results in stability of dc link voltage (i.e., by \eqref{Eq_theta_match}  $v_{\text{\upshape{dc}}}=\omega_{\text{GFC}}/{k_\theta}=\omega_{\text{SM}}/{k_\theta}$). The matching controlled converter switches its behavior as soon as $i_{\tau}$ exceeds the limit at approximately $t=0.5 \mathrm{s}$ in Figure \ref{fig:frequency_matching}. At around $t=3.5 \mathrm{s}$, the machine output power is sufficiently close to its steady-state value, $i_{{\text{\upshape{dc}}}}^\star$ and $i_{\tau}$ return to below the limit $i^{\text{dc}}_{\max}$, and the matching controlled converter recovers its dc voltage and frequency regulation capability and grid-forming dynamics. This behavior of matching control has been observed also for larger disturbance magnitudes. The nature of matching control - which accounts for the dc-side dynamics while regulating the ac dynamics - results in increased robustness with respect to large disturbances. In contrast, droop control, dVOC, and the VSM implicitly assume that the dc and ac sides are two independent systems and that can be regulated independently. This assumption is only justified under benign conditions and does not hold for large disturbances. As a consequence droop control, dVOC, and the VSM all exceed the limitations of the dc source for large disturbances and become unstable. 

\blue{Note that the behavior for the matching controlled GFCs is similar to that of the SM (see Section \ref{subsec:match} and \cite{cvetkovic_modeling_2015,arghir_grid-forming_2017,Catalin2019}), i.e., it achieves synchronization both under controlled or constant mechanical input power (i.e., dc current injection). While the VSM also enjoys some degrees of structural similarity to the SM (though only when the ac and dc currents are not saturated), its frequency dynamics are not linked to any physical storage element. In contrast, in the case of the matching control the dc voltage dynamics induce the frequency dynamics.} 

We observe the same instability of droop control, VSM, and dVOC when the test system contains one GFC and two SMs, i.e., the instability cannot be prevented by adding more inertia to the system. Figure \ref{fig:stable_idc_all_gfc} shows the dc current demand $i_\tau$ (i.e., before saturation) and dc voltage in an all-GFC system for a load increase of $\Delta p=0.9$ pu. The GFCs quickly synchronize to the post-event steady state, which does not exceed the maximum dc current, saturate the dc source for only approximately $200 \mathrm{ms}$, and remain stable. In contrast, in the system with two GFCs and one SM, the SM does not reach its increased post-event steady-state power injection for several seconds. During this time the response of droop control, VSM, and dVOC results in a power injection that exceeds the limits of the dc source and collapses the dc voltage. This highlights that the interaction of the fast GFC dynamics and slow SM dynamics contributes to the instability shown in Figure \ref{fig:unstable_idc}. \blue{Therefore, while matching control with dc-side limitation automatically limits the ac current in this particular scenario, the other grid-forming controls require additional mechanisms to limit the ac current.
}

\blue{\subsection{Impact of AC Current Limitation on GFCs}
In this subsection, we investigate if the ac current limitation \eqref{eq:ac_limiter} presented in section \ref{subsec:innercontr} can mitigate the instabilities observed in the previous section. To this end, we consider the same base load and disturbance as in previous test case, i.e., $p_l=2.25$ and $\Delta p=0.9$ pu. For this scenario, the GFCs dc transient current demand $i_\tau$ and the switching node current magnitude $\norm{\bm{i}_{\text{s,dq}}}$ exceeds the limits (i.e., $i^\text{dc}_{\max}=i^\text{ac}_{\max}=1.2$ pu) imposed by \eqref{eq:source_sat} and \eqref{eq:ac_limiter}. We remark that the combination of matching control and \eqref{eq:ac_limiter} is only included for the sake of completeness as it implicitly limits the ac current when dc source is saturated. 

We observe that using the ac current limiter does not stabilize the dc voltage of the GFCs controlled by droop control, VSM, and dVOC. Figure \ref{fig:dc and ac limiter} illustrates the behavior of the GFC at node 2 in Figure \ref{9bus_system}. Specifically, the current limitation imposed by \eqref{eq:ac_limiter} results in integrator windup, a loss of ac voltage control and, ultimately, instability of the grid-forming control \cite{Linbin} and dc voltage. We remark that GFCs exhibit the same instability behavior when the ac current limit is smaller than that of the dc-side i.e., $i^\text{ac}_{\max}<i^\text{dc}_{\max}$.  
%

To mitigate this load-induced instability, we explore a current limitation scheme that 
modifies the active power set-point when $\norm{\bm{i}_{\text{s,dq}}}$ exceeds a certain threshold value, i.e.,
\begin{align}\label{eq:dp}
\Delta p^\star\coloneqq\begin{cases}
0~&\If~\norm{\bm{i}_{\text{\upshape{s,dq}}}}\leq i^{\text{ac}}_\text{th},\\    
\gamma_p\left(\norm{\bm{i}_{\text{s,dq}}}-i^\text{ac}_\text{th}\right)~&\If~\norm{\bm{i}_{\text{\upshape{s,dq}}}}> i^{\text{ac}}_\text{th},
\end{cases}
\end{align}
and $\Delta p^\star$ is added to $p^\star$ in \eqref{Eq_idc_m}, \eqref{Eqtheta}, \eqref{eq:omega_vsm} and \eqref{Eq_dvoc_v}, $\gamma_p$ denotes a proportional control gain, and $i^\text{ac}_\text{th}<i^\text{ac}_{\max}$ is the activation threshold. Note that the control law \eqref{eq:dp} implicitly manipulates the grid-forming dynamics through their set-points such that the ac current magnitude stays within the admissible limits for large increases in load. We emphasize that this strategy aims at mitigating instabilities induced by large load increases and that the resulting GFC response to grid faults needs to be carefully studied.

For a $0.9$ pu load increase, the current limitation strategy \eqref{eq:dp} is able to stabilize the system with $i^{\text{ac}}_\text{th}=0.9$ pu and $\gamma_p=2.3 \left(p_\text{b}/i^\text{ac}_\text{b}\right)$ where $p_\text{b}$ and $i^\text{ac}_\text{b}$ denote the converter base power and current, respectively. Figure \ref{fig:dp} depicts the response of the same GFC as in Figure \ref{fig:dc and ac limiter}. Note that \eqref{eq:dp} effectively stabilizes the dc voltage for droop control, VSM, and dVOC. Moreover, in contrast to Figure \ref{fig:unstable_idc}, after the post-disturbance transient the dc source is no longer in saturation. 

Broadly speaking, this strategy succeeds to stabilize the system by steering the GFC power injection away from the critical limits. However, this also influences the post-disturbance operating point of the GFCs due to the threshold value being below the rated value. Thus, the GFCs no longer exhibit equal load-sharing with the SM. To further illustrate this Figure \ref{fig:mathcing vs VSM} depicts the power injection of VSMs with \eqref{eq:dp} and matching controlled GFCs without \eqref{eq:dp} (where each VSM is expected to inject $1.05$ pu in steady state). Figure \ref{fig:mathcing vs VSM} shows that the SM provides more power if  \eqref{eq:dp} is used for the VSMs. In contrast, matching controls inherent current limitation mechanism preserves the load-sharing capability. 
\begin{figure}[t!!]
    \centering
    \includegraphics[clip, trim=0.65cm 0.65cm 0.65cm 0.5cm,width=0.49\textwidth]{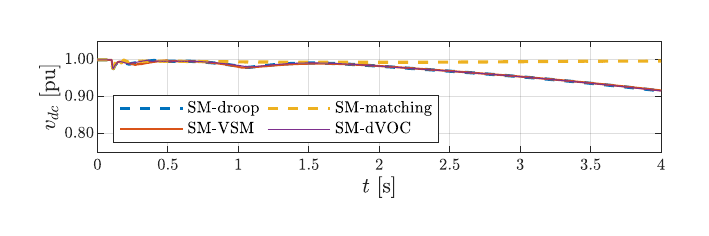}
    \caption{dc voltage of the converter at node 2 after a $0.9$ pu load disturbance when both dc and ac limitation schemes \eqref{eq:source_sat} and \eqref{eq:ac_limiter} are active and $\tau_\text{g}=5\mathrm{s}$.}
    \label{fig:dc and ac limiter}
\end{figure}
\begin{figure}[t!!]
    \centering
    \includegraphics[clip, trim=0.65cm 0.65cm 0.65cm 0.5cm,width=0.49\textwidth]{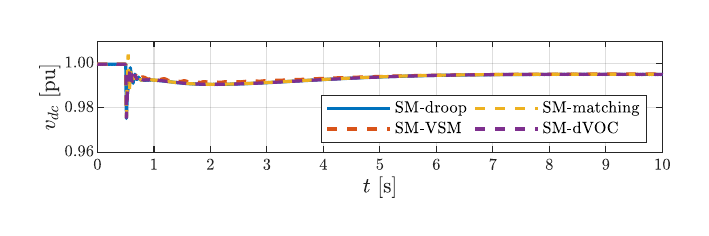}
    \includegraphics[clip, trim=0.65cm 0.65cm 0.65cm 0.5cm,width=0.49\textwidth]{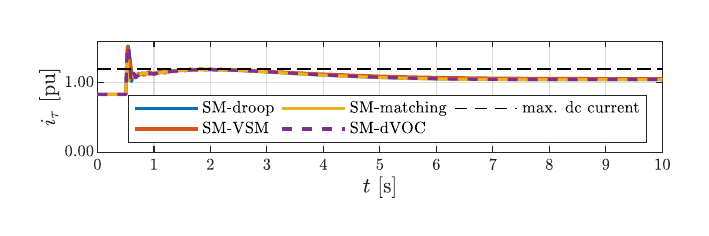}
    \includegraphics[clip, trim=0.65cm 0.65cm 0.65cm 0.5cm,width=0.49\textwidth]{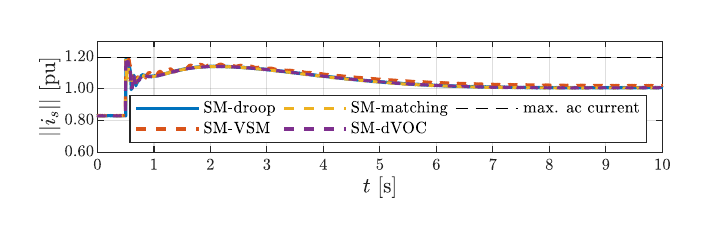}
    \caption{dc voltage (top left), dc current demand (top right) and ac current magnitude (bottom) of the converter at node 2 after a $0.9$ pu load disturbance when all the limitation schemes i.e., \eqref{eq:source_sat}, \eqref{eq:ac_limiter} and \eqref{eq:dp} are active and $\tau_\text{g}=5\mathrm{s}$.}
    \label{fig:dp}
\end{figure}
\begin{figure}[t!!]
    \centering
    \includegraphics[clip, trim=0.65cm 0.65cm 0.65cm 0.5cm,width=0.49\textwidth]{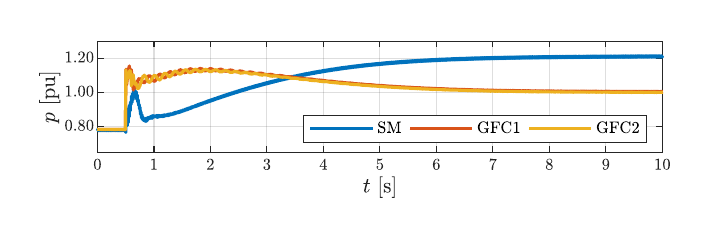}
    \includegraphics[clip, trim=0.65cm 0.65cm 0.65cm 0.5cm,width=0.49\textwidth]{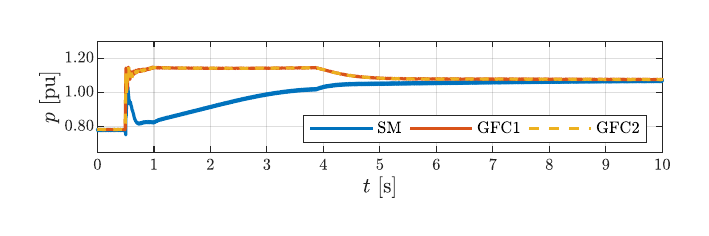}
    \caption{Active power injection of the SM and GFCs controlled by VSM technique with the control \eqref{eq:dp} (left), SM and GFCs controlled by matching technique without the control \eqref{eq:dp} (right) after a $0.9$ pu load disturbance.}
    \label{fig:mathcing vs VSM}
\end{figure}

Finally, we observe that the different time-scales in a low-inertia system contribute to the instabilities observed in this section. In particular, if the SM's turbine responds faster, GFCs with the standard limitation strategy \eqref{eq:ac_limiter} preserve stability - without the need to implement \eqref{eq:dp} - despite the fact that transient dc and ac currents exceed the limits. Figure \ref{fig:ac_lim_stable} shows the GFCs responses when the SM turbine delay $\tau_\text{g}$ is $1\mathrm{s}$ (cf. Figure \ref{fig:dc and ac limiter} where $\tau_\text{g}=5\mathrm{s}$). It can be seen that the slow SM turbine dynamics again contribute to the system instability when dc and ac currents are saturated. 
\begin{figure}[b!!]
    \centering
    \includegraphics[clip, trim=0.65cm 0.65cm 0.65cm 0.5cm,width=0.49\textwidth]{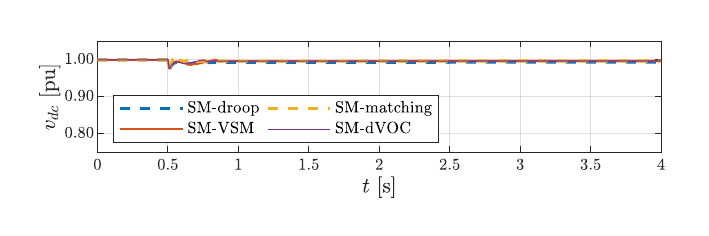}
    \caption{dc voltage of the converter at node 2 after a $0.9$ pu load disturbance when both dc and ac limitation schemes \eqref{eq:source_sat} and \eqref{eq:ac_limiter} are active and $\tau_\text{g}=1\mathrm{s}$.}
    \label{fig:ac_lim_stable}
\end{figure}

We conclude that the presence of different time-scales in a low-inertia system - often neglected in the literature \cite{WLBL15,denis_migrate_2018,ZK17,GDS15,Linbin,SHGM17,PD15,TWDB19} - must be considered in designing a robust ac current limitation mechanism for the GFCs. 
%
}

\subsection{Loss of Synchronous Machine Scenario} 
In this section, we study the response of grid-forming converters when disconnecting the synchronous machine at node 1, that is, the system turns into an all-GFCs network. The implications of such a contingency are threefold. First, the power injected by the machine, which partially supplies the base load, is no longer available. Second, the stabilizing dynamics associated with the machine's governor, AVR, and PSS are removed from the system. Third, the slow dynamics of the SM no longer interact with the fast dynamics of the GFCs. 

For this test, we set the base load to 2.1 pu, and the turbine and converter power set-points are set to $0.6$ and $0.75$ pu respectively. Note that when the SM at node 1 is disconnected, the converters increase their power output according to the power sharing behavior inherent to all four grid forming controls. The resulting increase in the converter power injection to roughly $1.05$ pu is similar to the load disturbance scenario used to study the unstable behavior of droop control, the VSM, and dVOC in the previous section.
\begin{figure}[t!!]
    \centering
    \includegraphics[clip, trim=0.65cm 0.65cm 0.65cm 0.5cm,width=0.49\textwidth]{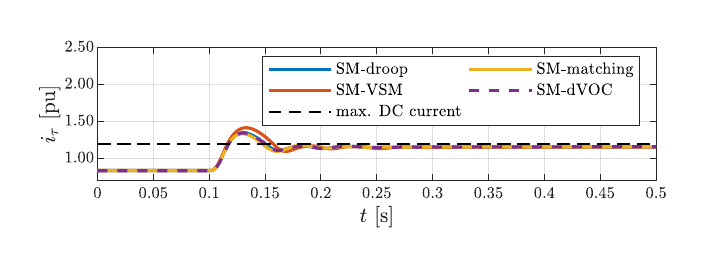}\hfill
    \includegraphics[clip, trim=0.65cm 0.65cm 0.65cm 0.5cm,width=0.49\textwidth]{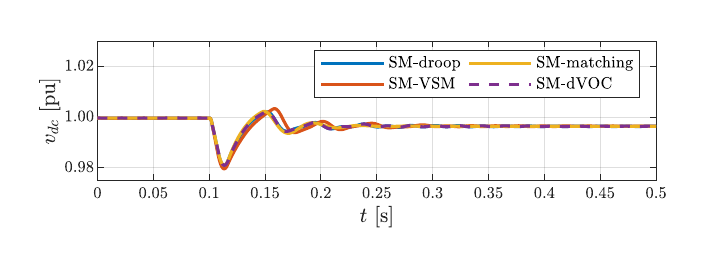}
    \caption{dc current demand ({\tb left}) and dc voltage ({\tb right}) of the converter at node 2 after the loss of the SM at node 1.}
    \label{fig:loss}
\end{figure}

Figure \ref{fig:loss} shows $i_{\tau}$ and $v_{{\text{\upshape{dc}}}}$ for the converter at node 2. Although the disturbance magnitude affecting the converters is similar to the one in studied in previous subsection, all GFCs remain stable after the loss of the SM. In particular, due to the absence of the slow turbine dynamics and fast synchronization of the converters $i_{\tau}$ is only above the limit $i^{\text{dc}}_{\max}$ for around $50 \mathrm{ms}$  while it remains above the limit for a prolonged period of time in Figure \ref{fig:unstable_idc}. 
This again highlights the problematic interaction between the fast response of the GFCs and the slow response of the SM. We stress that this adverse interaction is not resolved by increasing the inertia in the system (see the discussion the in previous subsections). While the synchronous machine perfectly meets classic power system control objectives on slower time scales, the dominant feature of GFCs is their fast response. However, the fast response of GFCs can also result in unforeseen interactions with other parts of the system such as the slow SM response (shown here), line dynamics (see \cite{VHH+18,gros_effect_2018}), and line limits \cite{aemo2019}. {For instance, during a recent separation event in Australia the rapid response of a battery energy storage system provided a valuable contribution to frequency stabilization but also contributed to tripping a line and islanding an area \cite[p. 67]{aemo2019}. Therefore, we expect to observe further adverse interactions in future studies related to low-inertia power system.}
\section{Qualitative analysis}\label{analysis}
In this section, we provide a qualitative but insightful analysis that explains the results observed in Section \ref{simulation_freq} and Section \ref{subsec:instability}. To this end, we develop simplified models that capture the small-signal frequency dynamics of synchronous machines and grid-forming converters. Applying arguments from singular perturbation theory \cite{PKC82,curi_control_2017} results in a model that highlights the main salient features of the interaction of synchronous machines and grid-forming converters.
\subsection{Frequency dynamics incorporating GFCs}
To obtain a simplified model of the frequency dynamics of the GFCs, we assume that the cascaded ac voltage and ac current control (see \eqref{Eq_vloop} and \eqref{Eq_iloop}) achieve perfect tracking (i.e., $\bm{v}_\text{\tb dq} = \hat{\bm{v}}_\text{\tb dq}$). Assuming that the system operates near the nominal steady-state (i.e., $\norm{\bm{v}_\text{\tb dq}}\approx v^\star$, $\omega \approx \omega^\star$) and $p^\star=0$, $q^\star=0$, we rewrite the remaining dynamics in terms of the voltage angle and power injection at every bus. This results in a simplified model of the angle $\underline{\theta}$ and the frequency $\underline{\omega}$ of a GFC or SM relative to a frame rotating at the nominal frequency $\omega^\star$.
\subsubsection{Droop control and dVOC}
For a converter controlled by droop control or dVOC, we obtain
\begin{align}\label{eq:reddroop}
 \dot{\underline{\theta}} &= - d_\omega p,
\end{align}
where $p$ is the power flowing out of the converter and $d_\omega$ is the droop control gain and given by $d_\omega=\eta / {v^\star}^2 $ for dVOC. 
\subsubsection{Synchronous machine, VSM, and matching control}\label{subsec:smmatch}
For a synchronous machine, a VSM, and a converter controlled by matching control we obtain
\begin{subequations}\label{eq:simpmatch}
\begin{align}
 \dot{\underline{\theta}} &= \underline{\omega}, \\
 2 H \dot{\underline{\omega}} &= -D \underline{\omega} + \sat\left(p_{\tau},p_{\max}\right)
 - p,\label{eq:simpmatch:freq}\\
   \tau \dot{p}_{\tau} &= -p_{\tau} - d_{p} \underline{\omega}.\label{eq:simpmatch:source}
\end{align}
\end{subequations}
For a SM the parameters directly correspond to the parameters of the machine model presented in Section \ref{syncmodel}, i.e., $H$, $d_p$, and $\tau=\tau_\text{g}$, are the machine inertia constant, governor gain, and turbine response time, and $D$ {\tb models the the effects of the damper windings \cite[Sec. 6.7]{SP98}}. Throughout this work we have not considered a limit on the turbine power output (i.e., $p_{\max} = \infty$) because a synchronous machine, in contrast to a GFC, typically has sufficient reserves to respond to the load changes and faults considered in this work. For the VSM presented in Section \ref{subsec:VSM}, we obtain $\tau = 0$, $p_{\tau}=0$, $H=1/2 J \omega^\star$, and $D=D_p \omega^\star$, i.e., the VSM does not emulate a turbine and implements no saturation of the damping term in its frequency dynamics \eqref{eq:omega_vsm}. Finally, for matching control we obtain $\tau = \tau_{\text{\upshape{dc}}}$, $p_{\max} = v^\star_{\text{\upshape{dc}}} i^{\text{dc}}_{\max}$, $H=1/2 C_{\text{\upshape{dc}}}/k^2_{\theta}$, {\tb $d_p=k_\text{dc}/k_\theta$} and $D=G_{\text{\upshape{dc}}}/k^2_{\theta}$ (see Section \ref{subsec:match}), i.e., by linking frequency and dc voltage, matching control clarifies that the dc source plays the role of the turbine in a machine and {\tb the proportional dc voltage control plays the role of a governor.}
\subsection{Reduced-order model}\label{subsec:redorder}
For brevity of the presentation we will now restrict our attention to the case of one SM and one GFC. The equivalent inertia constants and turbine time constants for the different grid-forming converter control strategies are either zero or negligible compared to typical inertia constants and turbine time constants for machines (see Table \ref{Table}). We therefore, assume that the states of the synchronous machine are slow variables, while the states of the GFC are fast and apply ideas from singular perturbation theory \cite{PKC82, curi_control_2017}. 

Using the dc power flow approximation and  $p_{d,\text{\upshape{GFC}}}$ and  $p_{d,\text{\upshape{SM}}}$ to denote a disturbance input, we obtain $p_{\text{\upshape{GFC}}} = b \left(\underline{\theta}_{\text{\upshape{GFC}}}-\underline{\theta}_{\text{\upshape{SM}}}\right) + p_{d,\text{\upshape{GFC}}}$ and $p_{\text{\upshape{SM}}} = b \left(\underline{\theta}_{\text{\upshape{SM}}}-\underline{\theta}_{\text{\upshape{GFC}}}\right) + p_{d,\text{\upshape{SM}}}$, where $b$ is the line susceptance. Neglecting the frequency dynamics and dc source dynamics, i.e., letting $\tau \to 0$ and $H \to 0$, the dynamics of the relative angle $\delta = \underline{\theta}_{\text{\upshape{SM}}}-\underline{\theta}_{\text{\upshape{GFC}}}$ are  given by $\dot{\delta} =  \underline{\omega}_{\text{\upshape{SM}}} - \left(D_{\text{\upshape{GFC}}}+d_{p,\text{\upshape{GFC}}}\right) \left(b \delta - p_{d,\text{\upshape{GFC}}}\right)$  if  $|d_{p,\text{\upshape{GFC}}}|<| p_{\max}|$, $\dot{\delta} =  \underline{\omega}_{\text{\upshape{SM}}} - D_{\text{\upshape{GFC}}} \left(b \delta - p_{d,\text{\upshape{GFC}}} \pm p_{\max}\right)$ if the dc source is saturated, and $D_{\text{\upshape{GFC}}}$ denotes is the damping provided by the GFC. For typical droop gains and network parameters, the relative angle dynamics are fast compared to the machine dynamics. Letting $\dot{\delta}\! \to \!0$, we obtain the reduced-order model
\begin{subequations}\label{eq:simpred}
\begin{align}
 \!\!2 H \dot{\underline{\omega}}_{\text{\upshape{SM}}} \!&=\! -D \underline{\omega}_{\text{\upshape{SM}}}\! -\! \sat\left(D_{\text{\upshape{GFC}}} \underline{\omega}_{\text{\upshape{SM}}} ,p_{\max}\right) \!+ \!p_{\tau,\text{\upshape{SM}}}
 \!+ \!p_d, \label{eq:simpred:freq}\\
   \!\tau \dot{p}_{\tau,\text{\upshape{SM}}} \!&=\! -p_{\tau,\text{\upshape{SM}}} - d_{p} \underline{\omega}_{\text{\upshape{SM}}}.\label{eq:simpred:turb}
\end{align}
\end{subequations}
where $p_d = p_{d,\text{\upshape{GFC}}} + p_{d,\text{\upshape{SM}}}$, $H$, $D$, $d_p$, and $\tau$ are the inertia constant, damping, governor gain, and turbine time constant of the synchronous machine, $D_{\text{\upshape{GFC}}}$ is the damping provided by the GFC, i.e., $D_{\text{\upshape{GFC}}} = 1/d_\omega$ (droop, dVOC, VSM) or  $D_{\text{\upshape{GFC}}} = d_p$ (matching). {Moreover, droop control, dVOC, and the VSM implement no saturation of their power injection in their respective angle / frequency dynamics (i.e., \eqref{eq:droop_omega}, \eqref{eq.droop.p.approx}, and \eqref{eq:omega_vsm}) resulting in $p_{\max}=\infty$. In contrast, for matching control the saturation of the dc source results in $p_{\max}=v^\star_{\text{\upshape{dc}}} i^{\text{dc}}_{\max}$.} 

We note that for the case of under damped dynamics \eqref{eq:simpred} and without saturation, a closed-form expression for the the step response and frequency nadir \eqref{eq:simpred} can be found in \cite[Sec. V-A]{PM18}. However, even for the seemingly simple model \eqref{eq:simpred} the dependence of the nadir on the parameters is very involved and does not provide much insight. By neglecting the damping term in \eqref{eq:simpred:freq} and the feedback term $-p_{\tau,\text{\upshape{SM}}}$ in \eqref{eq:simpred:turb} an insightful expression for the nadir is obtained in \cite{SER13}. However, the key feature of the GFCs is that they contribute damping, which is not captured by analysis in \cite{SER13}. Nonetheless, the model \eqref{eq:simpred} provides several insights that we discuss in the next section.
\subsection{Impact of Grid-Forming Control on Frequency Stability}
A system with three synchronous machines as in Section \ref{simulation} can be modeled by \eqref{eq:simpred} with $H=3 H_{\text{\upshape{SM}}}$, $d_p = 3 d_{p,\text{\upshape{SM}}}$, $\tau = \tau_{\text{\upshape{SM}}}$ and $D_{\text{\upshape{GFC}}}=0$. This corresponds to the well known \emph{center of inertia} frequency model with first-order turbine dynamics (see \cite{kundur1994power,PM18}). In contrast, if two SMs are replaced by GFCs with equal droop setting we obtain $H=H_{\text{\upshape{SM}}}$, $d_p = d_{p,\text{\upshape{SM}}}$, $\tau = \tau_{\text{\upshape{SM}}}$ and $D_{\text{\upshape{GFC}}}=2 d_{p,\text{\upshape{SM}}}$. 

This highlights that, on the time-scales of the SM, the GFCs provide fast acting frequency control. After an increase in load the machine inertia serves as buffer until the relatively slow turbine provides additional power to the machine. In contrast, the converters respond nearly instantaneously to any imbalance and therefore the need for inertia is decreased. Intuitively, this increase in fast primary frequency control should results in lower nadir values. Similarly, the additional damping provided by the converter in \eqref{eq:simpred} can be interpreted as a filter acting on the power imbalance, i.e., the GFCs are providing $D_{\text{\upshape{GFC}}} \underline{\omega}_{\text{\upshape{SM}}}$ or $p_{\max}$ and the power imbalance affecting the machine is reduced, therefore resulting in smaller average RoCoF. 

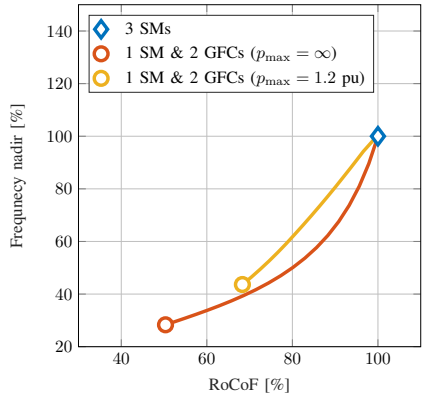
\begin{figure}[b!!!]
\centering
\newlength\fheight
\newlength\fwidth
\setlength{\fheight}{7cm}
\setlength{\fwidth}{7cm}
\input{comparison.tex}  
\centering
\caption{Change in averaged RoCoF and frequency nadir when transitioning from a system with 3 SMs to a system with one SM and two GFCs. \label{fig:comp}}    
\end{figure}

To validate the model \eqref{eq:simpred} and our interpretation, we compute the frequency nadir and averaged RoCoF (see \eqref{fmetrics}) for the machine parameters and disturbance used in in Section \ref{simulation}, and $H=\nu 3 H_{\text{\upshape{SM}}}$, $d_p = \nu 3 d_{p,\text{\upshape{SM}}}$, $\tau = \tau_{\text{\upshape{SM}}}$ and $D_{\text{\upshape{GFC}}}= (1-\nu) d_{p,\text{\upshape{SM}}}$ and $D=3 d_{p,\text{\upshape{SM}}}/10$, where $\nu \in [1/3,1]$ is a scalar parameter that interpolates the parameters between the two cases (all SM, one SM and two GFCs). The average RoCoF and frequency nadir according to \eqref{eq:simpred} are shown in Figure \ref{fig:comp}. The case with $p_{\max}=1.2$ corresponds to matching control, the one with $p_{\max}=\infty$ to droop control, dVOC, and VSMs. It can be seen that the GFCs result in an improvement compared to the all SM scenario, that the implicit saturation of the power injection by matching control results in a smaller improvement compared to droop control, dVOC and VSMs, and that the reduction in the average RoCoF and frequency nadir is line with the corresponding results in Figure \ref{all_rcf} and Figure \ref{all_mfd}. 

\subsection{Instability in the presence of large load disturbance}
The instabilities of droop control, dVOC, and the VSM observed in Section \ref{subsec:instability} can qualitatively be investigated using a simplified model of the dc-side. To compute the power $v_{\text{\upshape{dc}}} i_{\text{\upshape{x}}}$ flowing out of the dc-link capacitor, we assume that the controlled converter output filter dynamics are fast and can be neglected (i.e., $\dot{\bm{i}}_{\tb\text{\upshape{s}},\alpha\beta}=0$, $\dot{\bm{v}}_{\tb\alpha\beta}=0$) and that the ac output filter losses are negligible. This results in $v_{\text{\upshape{dc}}} i_{\text{\upshape{x}}} = \bm{v}^\top_{\tb\text{\upshape{s}},\alpha\beta} \bm{i}_{\tb\text{\upshape{s}},\alpha\beta} = \bm{v}^\top_{\tb\alpha\beta} \bm{i}_{\tb\alpha\beta} = p$. Moreover, we neglect the dc source dynamics to obtain the simplified dc voltage dynamics
\begin{align}
 C_{\text{\upshape{dc}}} \dot{v}_{\text{\upshape{dc}}} &= -G_{\text{\upshape{dc}}} v_{\text{\upshape{dc}}} + {\sat}\left(k_{\text{\upshape{dc}}}\left(v_{{\text{\upshape{dc}}}}^\star-v_{\text{\upshape{dc}}}\right),i^{\text{dc}}_{\max}\right)  - \dfrac{p}{v_{\text{\upshape{dc}}}},\label{eq:dcvoltred:dc}
\end{align} 
i.e., the active power $p$ flowing into the grid is drawn from the dc-link capacitor that is stabilized by a proportional control (see \eqref{Eq_idc_m}) if the dc current $i_{\tau}$ is not saturated. For a large enough constant perturbation $p>0$ the dc current in \eqref{eq:dcvoltred:dc} becomes saturated and controllability of the voltage $v_{\text{\upshape{dc}}}$ is lost and the dc voltage becomes unstable. 

In other words, if the dc source is saturated the power $p$ has to be controlled to stabilize the dc voltage. Matching control achieves this through the angle dynamics $\dot{\theta} = k_\theta v_{\text{\upshape{dc}}}$ which converge to a constant angle difference (i.e., $\underline{\omega}_{\text{\upshape{SM}}} = \underline{\omega}_{\text{\upshape{GFC}}}$) and power injection when the dc source is saturated and the SM has enough reserves to maintain $\underline{\omega}_{\text{\upshape{SM}}} \approx \omega^\star$ with the GFC providing its maximum output power. Moreover, for matching control it can be verified that $\left(\omega_{\text{\upshape{GFC}}} - \omega^\star\right)/\omega^\star = \left(v_{\text{\upshape{dc}}}-v_{{\text{\upshape{dc}}}}^\star\right)/v_{{\text{\upshape{dc}}}}^\star$. Therefore, the dc voltage deviation is proportional to the frequency deviation and the GFC with matching control remains stable for the scenario shown in Section \ref{subsec:instability}.

\section{Summary and Further Work}\label{conclusion}
In this paper we provided an extensive review of different grid-forming control techniques. Subsequently, we used the IEEE 9-bus test system incorporating high-fidelity GFC and SM models to compare the performance of different control techniques and their interaction with SM. Our case studies revealed that 1) integrating GFCs improves the frequency stability metrics compared to the baseline all-SMs system, 2) under a sufficiently large load disturbance some control techniques become unstable when the dc source current is limited for a prolonged time, 3) in contrast, matching control exhibits a change in operating mode from grid-forming to constant current source and remains stable, \blue{4) we explored the impact of the ac current limitation and proposed a solution to stabilize the grid-forming controls}, and 5)  we investigated the behavior of GFCs in response to the loss of a SM and highlighted a potentially destabilizing interaction between the fast synchronization of GFCs and the slow response of SMs. Moreover, we provided a qualitative analysis of the aforementioned results. Topics for our future works include investigating 1) further exploration of the impact of ac current limitation scheme for the GFCs, 2) the seamless transition between grid-forming and grid-following operation, 3) further theoretical analysis of the results observed in the case studies and 4) blending of the different control strategies into a controller that achieves their complementary benefits.
\appendices
\section{Tuning Criteria}\label{app:tunning}
The load-sharing capability of the control techniques presented in Section \ref{GFCC} is investigated in \cite{tayyebi_grid-forming_2018,darco_virtual_2013,colombino_global_2017}. Considering a heterogeneous network consisting of several GFCs (with different control) and SMs, we tune the control parameters such that all the units exhibit identical proportional load-sharing in steady-state. For the SM and droop controlled GFC, \eqref{Eq_sdroop} and \eqref{eq:droop_omega} can be rearranged to
\begin{subequations}\label{eq:ss_sm_and_droop}
    \begin{align}
    \omega^\star-\omega&=\frac{1}{d_p}\left(p-p^\star\right),\\
    \omega^\star-\omega&=d_\omega\left(p-p^\star\right).
    \end{align}
\end{subequations}
For VSM, assuming steady-state frequency and setting $\ddot{\theta}=\dot{\omega}=0$ in \eqref{eq:omega_vsm} results in
\begin{equation}\label{eq:ss_vsm}
\omega^\star-\omega=\dfrac{1}{D_p\omega^\star}\left(p-p^\star\right).
\end{equation}
For matching control, we assume that in steady-state $i_\text{dc} \approx i_\text{dc}^\star$ and $v_\text{dc}/v_\text{dc}^\star \approx1$. Setting $\dot{\omega}=0$ in \eqref{eq:vdc_matching} and replacing $i_\text{dc}$ by the expression from \eqref{Eq_idc_m} yields
\begin{equation}\label{eq:ss_matching}
\omega^\star-\omega=\frac{k_\theta}{k_\text{dc}v_\text{dc}}\left(p-p^\star\right).
\end{equation}
Lastly for dVOC, assuming $\norm{\hat{\bm{v}}_\text{\tb dq}} \approx v^\star$ in steady-state, the angle dynamics \eqref{eq.droop.p.approx} becomes
\begin{equation}\label{eq:ss_dvoc}
\omega^\star-\omega=\dfrac{\eta}{{v^\star}^2}\left(p-p^\star\right).
\end{equation}  
Hence, for any given droop gain $d_p$, if $d_\omega,D_p,k_\text{dc}~\text{and}~ \eta$ are selected such that the slopes of \eqref{eq:ss_sm_and_droop}-\eqref{eq:ss_dvoc} are equal, all the GFC control techniques and SM perform equal-load sharing. Moreover, by selecting $k_\text{dc}$ based on this criteria, the dc voltage control gain in \eqref{Eq_idc_m} is automatically set which is identical for all GFC implementation \cite{model}. 

Regarding the ac voltage regulation, the control gains in \eqref{Eq_vD}, \eqref{EQ_IF}, \eqref{EQ_VM}, and \eqref{eq.droop.q.approx} are selected to regulate the ac voltage at approximately equal time-scales. We refer to \cite{Prevost2018} for details on tuning the cascaded inner loops presented in Subsection \ref{subsec:innercontr}. It is noteworthy, that the time-scale of the reference model (i.e., grid-forming dynamics) must be slower than ac voltage control shown in Figure \ref{fig:architecture} to ensure optimal performance. Similarly, the ac current control must be faster than the outer voltage controller. Lastly, the choice of virtual inertia constant in \eqref{eq:omega_vsm} can largely influence VSM's dynamic behavior. We adopted the recommendation $J/D_p=0.02$ proposed in \cite{zhong_synchronverters:_2011}. The parameters used in our implementation \cite{model} are reported in Table \ref{Table}.           
{\fontfamily{ptm}\selectfont
    \begin{table}[h!]
        \centering
        \caption{Case study model and control parameters \cite{model}.\label{Table}}
        \begin{minipage}[c]{1\textwidth}
        \centering
        \scalebox{1.1}{
            {\renewcommand{\arraystretch}{1}    
                \begin{tabular}[]{|c|c||c|c||c|c|}
                    \hline\hline
                    \rowcolor{light-gray}
                    \multicolumn{6}{|c|}{{IEEE 9-bus test system base values}}\\
                    \hline
                    $S_\text{b}$ & $100$ MVA & $v_\text{b}$ & $230$ kV& $\omega_\text{b}$ & $2\pi50$ rad/s\\
                    \hline
                    \rowcolor{light-gray}
                    \multicolumn{6}{|c|}{{MV/HV transformer}}
                    \\ \hline
                    $S_\text{r}$ & $210$ MVA & $v_1$ & $13.8$ kV& $v_2$ & $230$ kV\\\hline
                    $R_1=R_2$ & $0.0027$ pu & $L_1=L_2$ & $0.08$ pu& $R_\text{m}=L_\text{m}$ & $500$ pu\\\hline
                    \rowcolor{light-gray}
                    \multicolumn{6}{|c|}{{single LV/MV transformer module in Figure \ref{agg_model}}}\\\hline
                    $S_\text{r}$ & $1.6$ MVA & $v_1$ & $1$ kV& $v_2$ & $13.8$ kV\\\hline
                    $R_1=R_2$ & $0.0073$ pu & $L_1=L_2$ & $0.018$ pu& $R_\text{m},L_\text{m}$ & $347,156$ pu\\\hline
                    \rowcolor{light-gray}
                    \multicolumn{6}{|c|}{{synchronous machine (SM)}}\\ \hline
                    $S_\text{r}$ & $100$ MVA & $v_\text{r}$ & $13.8$ kV& $\omega_\text{b}$ & $2\pi50$ rad/s\\\hline
                    $H$ & $3.7$ s & $d_p$ & $1$ \% & $\tau_{\text{g}}$ & $5$ s\\\hline
                    \rowcolor{light-gray}
                    \multicolumn{6}{|c|}{{single converter module in Figure \ref{agg_model}}}\\\hline
                    $S_\text{r}$ & $500$ kVA & $G_{{\text{\upshape{dc}}}},C_{{\text{\upshape{dc}}}}$ & $0.83,0.008$ $\Omega^{-1}$,F & $v^\star_{{\text{\upshape{dc}}}},v_\text{ll-rms}^\star$ & $2.44,1$ kV\\ \hline
                    $R$ & $0.001$ $\Omega$ & $L$ & $200$ $\mu$H & $C$ & $300$ $\mu$F\\\hline
                    $n$ & $100$ & $\tau_{{\text{\upshape{dc}}}}$ & $50$ ms & $i^{\text{dc}}_{\max}$ & $1.2$ pu\\\hline
                    \rowcolor{light-gray}
                    \multicolumn{6}{|c|}{{ac current, ac voltage, and dc voltage control}}\\\hline
                    $k_{v,\text{p}},k_{v,\text{i}}$ & $0.52,232.2$  &$k_{i,\text{p}},k_{i,\text{i}}$ & $0.73,0.0059$ & $k_{{\text{\upshape{dc}}}}$ & $1.6\times10^{3}$\\\hline
                    
                    \rowcolor{light-gray}
                    \multicolumn{6}{|c|}{{droop control}}\\\hline
                    $d_\omega$ & $2\pi0.05$ rad/s & $\omega^\star$ & $2\pi50$ & $k_\text{p},k_\text{i}$ & $0.001,0.5$\\\hline
                    \rowcolor{light-gray}
                    \multicolumn{6}{|c|}{{virtual synchronous machine (VSM)}}\\\hline
                    $D_p$ & $10^5$ & $J$ & $2\times10^3$  & $k_\text{p},k_\text{i}$ & $0.001,0.0021$ \\\hline
                    \rowcolor{light-gray}
                    \multicolumn{6}{|c|}{{matching control}}\\\hline
                    $k_\theta$ & $0.12$ &$k_{{\text{\upshape{dc}}}}$ & $1.6\times10^{3}$ & $k_\text{p},k_\text{i}$ & $0.001,0.5$ \\\hline
                    \rowcolor{light-gray}
                    \multicolumn{6}{|c|}{{dispatchable virtual oscillator control (dVOC)}}\\\hline
                    $\eta$ & $0.021$ & $\alpha$ & $6.66\times10^{4}$ & $\kappa$ & $\pi/2$  \\\hline
                    \hline
        \end{tabular}}}
        
        \end{minipage}
        \hfill
%
\end{table}}

\bibliographystyle{IEEEtran}
\bibliography{IEEEabrv,Ref}




\end{document}

%% file: comparison.tex
%
%

\definecolor{mycolor1}{rgb}{0.00000,0.44700,0.74100}%
\definecolor{mycolor2}{rgb}{0.85000,0.32500,0.09800}%
\definecolor{mycolor3}{rgb}{0.92900,0.69400,0.12500}%
\definecolor{Lgray}{gray}{0.75}

\begin{tikzpicture}[scale=0.65]

\begin{axis}[%
width=\fwidth,
height=\fheight,
at={(0\fwidth,0\fheight)},
scale only axis,
xmin=30,
xmax=110,
xlabel style={font=\color{white!15!black}},
xlabel={RoCoF [\%]},
ymin=20,
ymax=150,
legend style={row sep=2pt},
ylabel style={font=\color{white!15!black}},
ylabel={Frequnecy nadir [\%]},
axis background/.style={fill=white},
xmajorgrids,
ymajorgrids,
legend style={at={(0.03,0.97)}, anchor=north west, legend cell align=left, align=left, draw=white!15!black}
]
\addplot [only marks, color=mycolor1, line width=2.0pt, draw=none, mark size=5pt, mark=diamond*, mark options={solid, fill=white}]
  table[row sep=crcr]{%
100	100\\
};
\addlegendentry{\phantom{a}  3 SMs}


\addplot [only marks, color=mycolor2, line width=2.0pt, draw=none, mark size=4pt, mark=*, mark options={solid, fill=white}]
  table[row sep=crcr]{%
50.341887400956	28.3226886686405\\  
};
\addlegendentry{\phantom{a}  1 SM \& 2 GFCs ($p_{\max}=\infty$)}

\addplot [only marks, color=mycolor3, line width=2.0pt, draw=none, mark size=4pt, mark=*, mark options={solid, fill=white}]
  table[row sep=crcr]{%
68.3145100272789	43.6341936364971\\ 
};
\addlegendentry{\phantom{a}  1 SM \& 2 GFCs ($p_{\max}=1.2$ pu)}

%

\addplot [color=mycolor3, line width=2.0pt, forget plot]
  table[row sep=crcr]{%
100	100\\
99.6694339191595	99.7512468097569\\
98.9809733769911	98.6074737581797\\
97.9498962970934	96.7224126898384\\
96.5895766793943	94.244309141184\\
95.2119108305199	91.324105129123\\
93.6603468209229	88.0683472419819\\
91.9584294228594	84.579936713301\\
90.1790952048717	80.9475034031344\\
88.3144728745678	77.2391782496295\\
86.3937783658538	73.5115886288039\\
84.4393225473527	69.8088982102082\\
82.4670117111678	66.1646797503665\\
80.4940857399894	62.604278392458\\
78.5433813546293	59.1467873736625\\
76.5620633610397	55.800240853786\\
74.5488270873581	52.5720851435026\\
72.4976333977241	49.4661614071319\\
70.4373436679406	46.4881657578022\\
68.3145100272789	43.6341936364971\\
};

\addplot [color=mycolor2, line width=2.0pt, forget plot]
  table[row sep=crcr]{%
50.341887400956	28.3226886686405\\
52.5679836903396	29.5400030507456\\
54.9499804149878	30.8565228432792\\
57.4756919625359	32.2845015518607\\
60.126621927128	33.8382187040114\\
62.8798179296269	35.5344102580173\\
65.7098194336578	37.3928098002325\\
68.5903940538798	39.4368714890966\\
71.4959182282929	41.6946084346654\\
74.4023709891545	44.1997791026055\\
77.2879742747526	46.9933960450848\\
80.1335420901432	50.1257418450332\\
82.9226073383115	53.6590076558562\\
85.6413897571585	57.6710395125122\\
88.278657955443	62.2603883617928\\
90.8255269941412	67.5534536050699\\
93.2752223057906	73.7148316397322\\
95.6228317968435	80.9622359141718\\
97.8650609015662	89.5892610928551\\
100	100\\
};
\end{axis}
\end{tikzpicture}%